\documentclass[oneside,a4paper]{article}
\usepackage{amssymb,amsmath,latexsym,graphicx,tikz,indentfirst,scalerel,hyperref}
\def\dj{d\kern-0.4em\char"16\kern-0.1em}
\def\Dj{\mbox{\raise0.3ex\hbox{-}\kern-0.4em D}}
\newtheorem{teo}{Theorem}%[chapter]
\newtheorem{lem}[teo]{Lemma}
\newtheorem{tvr}[teo]{Proposition}
\newtheorem{pos}[teo]{Corollary}
\newtheorem{de}[teo]{Definition}
\newtheorem{pr}[teo]{Example}

\makeatletter
\renewcommand{\@dotsep}{10000}
\makeatother

\newenvironment{dokaz}
{\noindent
{\it Proof.}}
{\hspace{\stretch{1}}%
$\Box$}

\newcounter{primer}[section]

\DeclareMathOperator{\Inn}{Inn}

\DeclareMathOperator{\Aut}{Aut}
\DeclareMathOperator{\Z}{Z}
\DeclareMathOperator*{\bigboxtimes}{\scalerel*{\boxtimes}{\sum}}

\makeatletter 
\tikzset{my loop/.style =  {to path={
  \pgfextra{}
  [looseness=4,min distance=5mm]
  \tikz@to@curve@path},font=\sffamily\small
  }}  
\makeatletter

\newcommand{\lb}{\langle}
\newcommand{\rb}{\rangle}
\newcommand{\sime}{\stackrel{e}{\sim}}
\newcommand{\simeqe}{\stackrel{e}{\simeq}}
\newcommand{\equive}{\stackrel{e}{\equiv}}
\makeatletter
\newcommand{\myitem}[1]{%
\item[#1]\protected@edef\@currentlabel{#1}%
}
\makeatother

\begin{document}

\thispagestyle{empty}
\begin{center}
\Large{Enhanced Power Graphs of Finite Groups}

\normalsize{S. Zahirovi\' c, I. Bo\v snjak, R. Madar\' asz}
\end{center}

\begin{abstract}
The enhanced power graph $\mathcal G_e(\mathbf G)$ of a group $\mathbf G$ is the graph with vertex set $G$ such that  two vertices $x$ and $y$ are adjacent if they are contained in a same cyclic subgroup. We prove that finite groups with isomorphic enhanced power graphs have isomorphic directed power graphs. We show that any isomorphism between power graphs of finite groups is an isomorhism between enhanced power graphs of these groups, and we find all finite groups $\mathbf G$ for which $\Aut(\mathcal G_e(\mathbf G)$ is abelian, all finite groups $\mathbf G$ with $\lvert\Aut(\mathcal G_e(\mathbf G)\rvert$ being prime power, and all finite groups $\mathbf G$ with $\lvert\Aut(\mathcal G_e(\mathbf G)\rvert$ being square free. Also we describe enhanced power graphs of finite abelian groups. Finally, we give a characterization of finite nilpotent groups whose enhanced power graphs are perfect, and we present a sufficient condition for a finite group to have weakly perfect enhanced power graph.
\end{abstract}

\section{Introduction}

The directed power graph of a semigroup was introduced by Kelarev and Quinn \cite{prvi kelarev i kvin} and \cite{drugi kelarev i kvin} as the simple digraph whose vertex set is the universe of the semigroup, and with $x\rightarrow y$ if $y$ is a power of $x$. The (undirected) power graph of a semigroup was defined by Chakrabarty in \cite{cakrabarti} as the simple graph whose vertex set is the universe of the semigroup, and with $x$ and $y$ being adjacent if one of them is a power of the other. In \cite{power graph} Cameron and Ghosh proved that finite groups with isomorphic directed power graphs have the same numbers of elements of each order, and that finite abelian groups with isomorphic power graphs are isomorphic. Later in \cite{power graph 2} Cameron proved that finite groups with isomorphic power graphs have isomorphic directed power graphs. The authors of \cite{on the structure} and \cite{feng ma i vang} independently proved that the power graph of any group of bounded exponent is perfect. In \cite{sitov} Shitov proved that chromatic number of the power graph of any semigroup is at most countable. Commuting graph of a group appears to be defined in \cite{brauer i foler}. It is the simple graph with the universe of the group being its vertex set, in which $x$ and $y$ are joined if $x$ and $y$ commute. Note that usually the vertex set of the commuting graphs contains only non-central elements of the group. In \cite{on the structure} the authors introduced the enhanced power graph of a group whose vertex set is the universe of the group, and with $x$ and $y$ being adjacent if $x$ and $y$ are contained in a same cyclic group. Note that the enhanced power graph of a group lies, as a subgraph, between the power graph and the commuting graph of the group, if the vertex set of the commuting graph is entire universe of the group. In \cite{on the structure} the authors proved that maximal cliques of the enhanced power graph of a group is a cyclic subgroup or a locally cyclic subgroup of the group.

In this paper we study properties of the enhanced power graph of a finite group, and its relations to the power graph and the directed power graph of the group.

In Section \ref{osnovne osobine} we provide some basic definitions and notations which will be used in this paper. Section \ref{power graf i enhanced power graf} is motivated by the results of \cite{power graph} and \cite{power graph 2}, and in this section we prove that finite groups with isomorphic enhanced power graphs have isomorphic directed power graphs. This implies that finite groups have isomorphic enhanced power graphs if and only if they have isomorphic directed power graphs  if and only if  they have isomorphic power graphs. We also give an efficient way to determine the number of elements of the group of any order from the enhanced power graph. Finally, we show that it is possible to reconstruct any finite abelian group from its enhanced power graph.

Section \ref{grupe automorfizama} is motivated by \cite{power graph}. Here we extend Cameron's proof from \cite{power graph 2} in order to show that an isomorphism of the power graph of a finite group is an isomorphism of the enhanced power graph of that group too. We also find all finite groups $\mathbf G$ with abelian $\Aut(\mathcal G_e(\mathbf G))$, all finite groups $\mathbf G$ such that $|\Aut(\mathcal G_e(\mathbf G))|$ is square free, and all finite groups $\mathbf G$ such that $|\Aut(\mathcal G_e(\mathbf G))|$ is prime power.

In Section \ref{p semistabla} we give a description of enhanced power graphs of finite abelian $p$-groups, which enables us to describe enhanced power graphs of all finite abelian groups. Given a finite graph, we provide some necessary conditions for it to be the enhanced power graph of a group. Also, we give an algorithm which for a finite graph satisfying these conditions, were it the enhanced power graph of a finite group with a unique $p$-Sylow subgroup, returns the enhanced power graph of its $p$-Sylow subgroup.

Finally in Section \ref{perfektnost}, which is motivated by \cite{on the structure}, shows that there are groups which do not have perfect enhanced power graph. Furthermore, we give a characterization of finite nilpotent groups whose enhanced power graphs are perfect. Also for a class of groups, which contains all finite nilpotent groups,  we prove that all of them have weakly perfect enhanced power graphs.

\section{Basic notions and notations}
\label{osnovne osobine}

{\bf Graph} $\Gamma$ is a structure $(V(\Gamma),E(\Gamma))$, or simply $(V,E)$, where $V$ is a set, and $E\subseteq V^{[2]}$ is a set of two-element subsets %\footnote{Umesto $V^{[2]}$ mogu jo\v s da se koriste i oznake $\binom{V}{2}$, ili $[V]^2$. $\binom{V}{2}$ ne bi bilo najzgodnije re\v senje, jer onda ne bi mo\v zda bilo jasno da li je $\binom nk$ skup $k$-to\v clanih podskupova, ili broj.}
  of $V$. Here $V$ is the set of vertices, and $E$ is the set of edges. For $x,y\in V$ we say that $x$ and $y$ are adjacent in $\Gamma$ if $\{x,y\}\in E$, and we denote it with $x\sim_\Gamma y$, or simply $x\sim y$. We say that graph $\Delta=(V_1,E_1)$ is a {\bf subgraph} of graph $\Gamma=(V_2,E_2)$ if $V_1\subseteq V_2$ and $E_1\subseteq E_2$. $\Delta$ is an {\bf induced subgraph} of $\Gamma$ if $V_1\subseteq V_2$, and $E_1=E_2\cap V_1^2$. In this case we also say that $\Delta$ is an induced subgraph of $\Gamma$ by $E_1$, and we denote that with $\Delta=\Gamma[E_1]$. {\bf Strong product} of graphs $\Gamma$ and $\Delta$ is the graph $\Gamma\boxtimes\Delta$ such that
\begin{align*}
(x_1,y_1)\sim_{\Gamma\boxtimes\Delta}(x_2,y_2)\text{ if }&(x_1=x_2\wedge y_1\sim_\Delta y_2)\\
&\vee(x_1\sim_\Gamma x_2\wedge y_1= y_2)\\
&\vee(x_1\sim_\Gamma x_2\wedge y_1\sim_\Delta y_2).
\end{align*}

{\bf Closed neighborhood} of a vertex $x$ in a graph $\Gamma$, denoted with $\overline N_\Gamma(x)$, or simply $\overline N(x)$, is the set $\overline N_\Gamma(x)=\{y\mid y\sim_\Gamma x\text{ or } y=x\}$.

{\bf Directed graph} (or {\bf digraph}) $\vec \Gamma$ is a structure $(V(\vec\Gamma),E(\vec\Gamma))$, or simply $(V,E)$, where $V$ is a set, and $E$ is an irreflexive relation on $V$. Here two $V$ and $E$ are the set of vertices and the set of edges, respectively. If $(x,y)\in E$, then we say that $y$ is a direct successor of $x$, or that $x$ is a direct predecessor of $y$, and we denote that with $x\rightarrow_{\vec\Gamma} y$, or simply $x\rightarrow y$.

Throughout this paper, we shall denote algebraic structures, such as groups, with bold capitals, while we will denote their universes with respective regular capital letters. For $x,y\in G$, where $\mathbf G$ is a group, we shall write $x\approx_{\mathbf G} y$, or simply $x\approx y$, if $\langle x\rangle=\langle y\rangle$. Also, $o(x)$ shall denote the order of an element $x$ of a group.

\begin{de}
{\bf Directed power graph} of a group $\mathbf G$ is the digraph $\vec{\mathcal G}(\mathbf G)$ whose vertex set is $G$, and  there is a directed edge from $x$ to $y$, $y\neq x$,  if $y\in\lb x\rb$. We shall denote that there is a directed edge from $x$ to $y$ in $\vec{\mathcal G}(\mathbf G)$ with $x\rightarrow_\mathbf G y$, or simply $x\rightarrow y$.

{\bf Power graph} of a group $\mathbf G$ is the graph $\mathcal G(\mathbf G)$  whose vertex set is $G$, and  $x$ and $y$, $x\neq y$, are adjacent if either $y\in\lb x\rb$ or $x\in\lb y\rb$. We shall denote that $x$ and $y$ are adjacent in $\mathcal G(\mathbf G)$ with $x\sim_\mathbf G y$, or simply $x\sim y$.

We shall write  $x\equiv_{\mathbf G} y$, or simply $x\equiv y$, if $\overline N_{\mathcal G(\mathbf G)}(x)=\overline N_{\mathcal G(\mathbf G)}(y)$.

{\bf Enhanced power graph} of a group $\mathbf G$ is the graph $\mathcal G_e(\mathbf G)$  whose vertex set is $G$, and such that $x$ and $y$, $x\neq y$, are adjacent if there exists $z\in G$ such that $x,y\in\lb z\rb$.  We shall denote that $x$ and $y$ are adjacent in $\mathcal G_e(\mathbf G)$ with $x\sime_\mathbf G y$, or simply $x\sime y$.%\footnote{Obratiti pa\v znju na ove definicije zbog grupa sa elementima beskona\v cnog reda. Kameron je u  radu sa torziono slobodnim grupama ve\' c na sli\v can nacin definisao grupe, ali je on izbacivao slu\v caj $x^0=y$, pa elementi beskona\v cnog reda nisu bili susedni sa jedini\v cnim elementom. Ipak, napominjao je u tom radu da bi isto sve radilo i sa definicijom kao ovom sad jer se jedini\v cni element (makar kod njega) uvek mo\v ze prona\' ci bar do na izomorfizam.}

We shall write $x\equive_{\mathbf G}y$, or simply  $x\equive y$, if $\overline N_{\mathcal G_e(\mathbf G)}(x)=\overline N_{\mathcal G_e(\mathbf G)}(y)$.
\end{de}

\begin{lem}\label{proizvodi}
Let $\mathbf G$ and $\mathbf H$ be groups. Then:
\begin{enumerate}
\item $\mathcal G_e(\mathbf G\times\mathbf H)\subseteq \mathcal G_e(\mathbf G)\boxtimes\mathcal G_e(\mathbf H)$;
\item Let $\mathbf G$ and $\mathbf H$ be torsion groups. Then $\mathcal G_e(\mathbf G\times\mathbf H)= \mathcal G_e(\mathbf G)\boxtimes\mathcal G_e(\mathbf H)$  if and only if  $\gcd(o(g),o(h))=1$ for all $g\in G$ and $h\in H$.
\end{enumerate}
\end{lem}

\begin{dokaz}
{\it 1.} If $(x_1,y_1)\sime_{\mathbf G\times\mathbf H}(x_2,y_2)$, then  $(x_1,y_1),(x_2,y_2)\in\lb (x_3,y_3)\rb$, which implies $x_1,x_2\in\lb x_3\rb$ and $y_1,y_2\in\lb y_3\rb$. Now $(x_1,y_1)\sim_{\mathcal G_e(\mathbf G)\boxtimes\mathcal G_e(\mathbf H)}(x_2,y_2)$ trivially follows, since we have $x_1\neq x_2$ or $y_1\neq y_2$. This proves the inclusion.

{\it 2.} Suppose that $\gcd(o(g),o(h))=1$ for all $g\in G$ and $h\in H$.  One inclusion has already been proved, so we only need to prove the other one. Let us assume that $(x_1,y_1)\sim_{\mathcal G_e(\mathbf G)\boxtimes\mathcal G_e(\mathbf H)}(x_2,y_2)$. This implies $x_1,x_2\in\lb x_3\rb$ and $y_1,y_2\in\lb y_3\rb$ for some $x_3\in  G$ and $y_3\in  H$. Since $x_3$ and $y_3$ have relatively prime orders, it follows $\lb (x_3,1),(1,y_3)\rb=\lb(x_3,y_3)\rb$, which implies $(x_1,y_1),(x_2,y_2)\in\lb(x_3,y_3)\rb$. Because either $x_1\neq x_2$ or $y_1\neq y_2$, we have $(x_1,y_1)\sime_{\mathbf G\times\mathbf H}(x_2,y_2)$, which proves the other inclusion. This proves one implication. Let us prove the other.

Assume that there are elements $g_0\in G$ and $h_0\in H$ with $\gcd(o(g_0),o(h_0))=k>1$. Since $\big(\mathcal G_e(\mathbf K)\big)[L]=\mathcal G_e(\mathbf L)$ for any groups $\mathbf K$ and $\mathbf L$ with $\mathbf L\leq\mathbf K$, we have:
\begin{align*}
\big(\mathcal G_e(\mathbf G\times\mathbf H)\big)[\lb g_0\rb\times \lb h_0\rb]=\mathcal G_e(\lb g_0\rb\times \lb h_0\rb)<K_{o(g_0)\cdot o(h_0)}=K_{o(g_0)}\boxtimes K_{o(h_0)}\\
=\mathcal G_e(\lb g_0\rb)\boxtimes\mathcal G_e(\lb h_0\rb)=\big(\mathcal G_e(\mathbf G)\boxtimes\mathcal G_e(\mathbf H)\big)[\lb g_0\rb\times \lb h_0\rb].
\end{align*}
Therefore, $\mathcal G_e(\mathbf G\times\mathbf H)= \mathcal G_e(\mathbf G)\boxtimes\mathcal G_e(\mathbf H)$, which proves the other implication as well.
\end{dokaz}\\

%We shall use the following notations. $x\rightarrow_\mathbf G y$ or $x\rightarrow y$ if x
%\begin{align*}
%&x\rightarrow y\text{ iff }y\text{ is a power of }x;\\
%&x\sim y\text{ iff }x\rightarrow y\text{ or }x\leftarrow y;\\
%&x\sime y\stackrel{e}{\sim}\text{ iff  there exists }z\text{ so that }z\rightarrow x\text{ and }z\rightarrow y.
%\end{align*}

\section{Enhanced power graph and the power graph}
\label{power graf i enhanced power graf}

It is obvious that, if two groups have isomorphic directed power graphs, then they have isomorphic power graphs. Peter Cameron proved in \cite{power graph 2} that the opposite holds too.

\begin{teo}\label{neusmereni odredjuju usmerene}
(Theorem 2 in \cite{power graph 2}). If $\mathbf G$ and $\mathbf H$ are finite groups whose power graphs are isomorphic, then their directed power graphs are also isomorphic.
\end{teo}

 We will prove that, if two finite groups have isomorphic enhanced power graphs, then they have isomorphic directed power graphs too. One can easily see that the opposite implication holds as well.

Note that there are non-isomorphic groups which have isomorphic power graphs. Peter Cameron and Shamik Ghosh gave several examples of such groups in \cite{power graph}, with some of them being finite. By our results from this section these finite groups  have isomorphic enhanced power graphs. Before proceeding we will introduce the following lemma which was proved in \cite{on the structure}.

\begin{lem}\label{klike i ciklicne grupe}
(Lemma 36 in \cite{on the structure}). A maximal clique in the enhanced power graph is either a cyclic subgroup or a locally cyclic subgroup.
\end{lem}

\begin{lem}\label{uredjene klase}
Let $\mathbf G$ be a finite group, and let $\mathcal Cl$ be the set of all maximal cliques of $\mathcal G_e(\mathbf G)$. Then for any $x,y\in G$:
\begin{enumerate}
\item $\overline N(x)\subseteq\overline N(y)$  in $\mathcal G_e(\mathbf G)$ if and only if  $x\in C\Rightarrow y\in C$ for all $C\in \mathcal Cl$;
\item  $x\equive y$  if and only if  $x\in C\Leftrightarrow y\in C$ for all $C\in \mathcal Cl$;
\item If $\langle x\rangle=\langle y\rangle$, then $x\equive y$;
\item Relation $\leq$ defined with: $[x]_{\equive}\leq[y]_{\equive}$  if and only if  $y\in C\Rightarrow x\in C$ for all $C\in \mathcal Cl$, is an order on $G/\equive$;
\item $\langle[x]_{\equive}\rangle=\bigcap\{C\in\mathcal Cl\mid [x]_{\equive_{\mathbf G}}\subseteq C\}$;
\item Any $\equive$-class contains either $\varphi(m)$ or no elements of order $m$.\footnote{$\varphi$ represents the Euler function.}

$[x]_{\equive}$ contains $\varphi(m)$ elements of order $m$  if and only if  it is a minimal $\equive$-class (with respect to $\leq$ introduced in {\it 4.}) with property that $m\big||\langle[x]_{\equive}\rangle|$.
\end{enumerate}
\end{lem}

\begin{dokaz}
{\it 1.} $(\Rightarrow)$ Let $x,y,z\in G$. Suppose that $\overline N(x)\subseteq\overline N(y)$, and $x\in C$ for some $C\in\mathcal Cl$. By Lemma \ref{klike i ciklicne grupe} $\mathbf C$ is a cyclic subgroup of $\mathbf G$. Let $\mathbf C=\langle z\rangle$. Then $y\sime z$ because $x\sime z$. But as $\langle z\rangle$ is a maximal cyclic subgroup, we have $z\rightarrow y$, which implies $y\in C$.

$(\Leftarrow)$ Let $x,y\in G$, and suppose that $x\in C\Rightarrow y\in C$ for all $C\in \mathcal Cl$. If $z\sime x$, then $z,x\in C$ for some maximal cyclic subgroup $\mathbf C$ of $\mathbf G$. Then by Lemma \ref{klike i ciklicne grupe} $C\in\mathcal Cl$, and therefore $z,y\in C$, which implies $z\sime x$.

{\it 2.} and {\it 4.} follow from {\it 1.}, while {\it 3.} obviously holds. 

{\it 5.} Let $x\in G$, and let us denote $\bigcap\{C\in\mathcal Cl\mid [x]_{\equive}\subseteq C\}$ with $D_x$. $\mathbf D_x$ is a cyclic subgroup, so let $\mathbf D_x=\langle z\rangle$. It is sufficient to prove that $x\equive_{\mathbf G}z$. $z\in D_x$ implies $x\in C\Rightarrow z\in C$ for all $C\in\mathcal Cl$. On the other hand $z\in C\Rightarrow x\in C$ for all $C\in\mathcal Cl$, because $x\in\langle z\rangle$. Therefore, $x\equive z$.

{\it 6.} If a $\equive$-class contains an element $z_0$ of order $m$, then it contains all $z$ with $\langle z\rangle=\langle z_0\rangle$. Also, by {\it 5.} $[z_0]_{\equive}$ doesn't contain any element $t$ of order $m$ such that $\langle z_0\rangle\neq\langle z\rangle$. For that reason $[z_0]_{\equive}$ contains exactly $\varphi(m)$ elements of order $m$.

Let $\mathcal D$ be the set of $\equive$-classes $[y]_{\equive_{\mathbf G}}$ with the property $m\big||\langle[y]_{\equive}\rangle|$.
%Let $[y]_\equiv$ be minimal in $\mathcal D_m$. By {\it 4.} $\bigcap\{C\in\mathcal Cl^G\mid [y]_{\equive_{\mathbf G}}\subseteq C\}=\langle[y]_{\equive_{\mathbf G}}\rangle$.
%
Suppose that $[x]_{\equive}$ is minimal in $\mathcal D$. By {\it 5.} there is an element $z\in\langle[x]_{\equive}\rangle$ with $o(z)=m$, and by {\it 4.} and {\it 5.} it follows $[z]_{\equive}\leq[x]_{\equive}$. Were it $[z]_{\equive}<[x]_{\equive}$, then $[x]_{\equive}$ would not be minimal in $\mathcal D$. Therefore, $z\in [x]_{\equive}$. Let us prove now the other implication as well. If $y\in[x]_{\equive}$ with $o(y)=m$, then $\langle y\rangle\leq\langle\mathbf [x]_{\equive}\rangle$, and $[x]_{\equive}\in\mathcal D$. Suppose that $[x]_{\equive}$ is not minimal in $\mathcal D$. Then there is $x_0\in G$ with $[x_0]_{\equive}<[x]_{\equive}$ and $[x_0]_{\equive}$ being minimal in $\mathcal D$. Therefore, there is $y_0\in[x_0]_{\equive}$ with $o(y_0)=m$, and we have $\langle y_0\rangle\leq\langle[x_0]_{\equive}\rangle\leq \langle[x]_{\equive}\rangle$. Also, since $y_0\not\equive y$, by {\it 3.} we get $\langle y_0\rangle\neq\langle y\rangle$, and $\langle y\rangle$ and $\langle y_0\rangle$ are two different cyclic subgroups of order $m$ of the cyclic subgroup  $\langle[x_0]_{\equive}\rangle$. This is a contradiction, so $[x]_{\equive}$ is minimal in $\mathcal D$, and the equivalence has been proved.
\end{dokaz}

\begin{teo}\label{enhanced odredjuje usmereni}
Let $\mathbf G$ and $\mathbf H$ be finite groups. If $\mathcal G_e(\mathbf G)\cong\mathcal G_e(\mathbf H)$, then $\vec{\mathcal G}(\mathbf G)\cong\vec{\mathcal G}(\mathbf H)$.
\end{teo}

\begin{dokaz}
Let $\mathbf G$ and $\mathbf H$ be finite groups with $\mathcal G_e(\mathbf G)\cong\mathcal G_e(\mathbf H)$, and let $\psi$ be the isomorphism. %Let $\mathcal Cl^G=\{C_0^G,C_1^G,...,C_{n-1}^G\}$ be all maximal cliques of $\mathcal G_e(\mathbf G)$. Also let $\mathcal Cl^H=\{C_0^H,C_1^H,...,C_{n-1}^H\}$ be the corresponding maximal cliques of $\mathcal G_e(\mathbf H)$ by the enhanced graph isomorphism. Then by Lemma \ref{klike i ciklicne grupe}, $\mathbf C_0^G, \mathbf C_1^G,...,\mathbf C_{n-1}^G$  are all maximal subgroups of $\mathbf G$, and $\mathbf C_0^H, \mathbf C_1^H,...,\mathbf C_{n-1}^H$  are all maximal subgroups of $\mathbf H$. For $X\subseteq G$ and $Y\subseteq H$ of the same sizes we will say $X$ and $Y$ are corresponding iff $X\subseteq C_i^G\Leftrightarrow Y\subseteq C_i^H$ for all $i\in n$. We shall say $x\equive_{\mathbf G} y$ iff $x$ and $y$ have the same closed neighborhoods in $\mathcal G_e(\mathbf G)$. In the same manner we define $\equive_{\mathbf H}$. Obviously $\equive_{\mathbf G}$ and $\equive_{\mathbf H}$ are equivalence relations. Notice that, for any $x,y\in G$, $x\equive_{\mathbf G} y$, iff $x\in C\Leftrightarrow y\in C$ for all $C\in\mathcal Cl^G$. We order $\equive_{\mathbf G}$-classes with $[x]_{\equive_{\mathbf G}}\leq[y]_{\equive_{\mathbf G}}$ iff  $y\in C\Rightarrow x\in C$ for all $C\in\mathcal Cl^G$, which is equivalent to $\overline N_{\mathcal G_e(\mathbf G)}(y)\subseteq\overline N_{\mathcal G_e(\mathbf G)}(x)$. We shall also say $x\approx_{\mathbf G}$ iff $\langle x\rangle=\langle y\rangle$, for any $x,y\in G$, and $\approx_{\mathbf H}$ is defined in the same manner. It is easy to see that $\approx_{\mathbf G}$ and $\approx_{\mathbf H}$ are refinements of $\equive_{\mathbf G}$ and $\equive_{\mathbf H}$, respectively.
Let $C_1^G,C_2^G,...,C_n^G$ be all maximal cliques of $\mathcal G_e(\mathbf G)$, and let $C_1^H,C_2^H,...,C_n^H$ be all maximal cliques of $\mathcal G_e(\mathbf H)$ with $\psi(C_i^G)=C_i^H$ for all $i\leq n$. For $X\subseteq G$ and $Y\subseteq H$ of the same sizes we will say $X$ and $Y$ are corresponding if $X\subseteq C_i^G\Leftrightarrow Y\subseteq C_i^H$ for all $i\leq n$. %For any $x,y\in G$ let $x\approx_{\mathbf G} y$ if $\langle x\rangle=\langle y\rangle$, and let $\approx_{\mathbf H}$ be defined in similar way. 

Let $[x]_{\equive_{\mathbf G}}$ be a $\equive_{\mathbf G}$-class,  $y=\psi(x)$, and let $m\in \mathbb N$. Then $[y]_{\equive_{\mathbf H}}$ is corresponding $\equive_{\mathbf H}$-class to $[x]_{\equive_{\mathbf G}}$. By Lemma \ref{uredjene klase}, for any $m\in \mathbb N$, $[x]_{\equive_{\mathbf G}}$ contains an element of order $m$  if and only if  $[y]_{\equive_{\mathbf H}}$ does, while both $[x]_{\equive_{\mathbf G}}$ and $[y]_{\equive_{\mathbf G}}$ contain at most one $\approx$-class made up of elements of order $m$. Therefore, for any $\approx_{\mathbf G}$-class $[x_0]_{\approx_{\mathbf G}}$ there is a corresponding $\approx_{\mathbf H}$-class $[y_0]_{\approx_{\mathbf H}}\subseteq [y_0]_{\equive_{\mathbf H}}$ with $o(x_0)=o(y_0)$. Let $\vartheta:G\rightarrow H$ be a bijection which maps any $\approx_{\mathbf G}$-class $[x]_{\approx_{\mathbf G}}$ onto its corresponding $\approx_{\mathbf H}$-class $[y]_{\approx_{\mathbf H}}$ with $o(x)=o(y)$. Let us prove that $\vartheta$ is a digraph isomorphism from $\vec{\mathcal G}(\mathbf G)$ to $\vec{\mathcal G}(\mathbf H)$. Suppose that $x_0\rightarrow_{\mathbf G}x_1$, and $\vartheta(x_0)=y_0$ and $\vartheta(x_1)=y_1$. Then $o(x_1)|o(x_0)$, and $x_0,x_1\in C_i^G$, for some $i\in n$. Since $y_0$ and $y_1$ are in corresponding $\approx_{\mathbf H}$-classes to $\approx_{\mathbf G}$-classes of $x_0$ and $x_1$, respectively, we get $y_0,y_1\in C_i^H$. Because $y_0$ and $y_1$ are contained in a cyclic subgroup, and $o(y_1)|o(y_0)$, it follows that $y_0\rightarrow_{\mathbf H}y_1$. This proves the theorem.
%
%
%By Lemma \ref{uredjene klase} we can decide from $\mathcal G_e(\mathbf H)$ whether a $\equive_{\mathbf H}$-class contains $0$ or $\varphi(m)$ elements of order $m$ as well. For $\equive_{\mathbf G}$-class $[x]_{\equive_{\mathbf G}}$ there is a corresponding $\equive_{\mathbf H}$-class $[y]_{\equive_{\mathbf H}}$, and by the above shown, for any $\approx_{\mathbf G}$-class $[x_0]_{\approx_{\mathbf G}}\subseteq[x]_{\equive_{\mathbf G}}$ there is a corresponding $\approx_{\mathbf H}$-class $[y_0]_{\approx_{\mathbf H}}\subseteq [y]_{\equive_{\mathbf H}}$ with $o(x_0)=o(y_0)$. Let $\vartheta:G\rightarrow H$ be a bijection which maps any $\approx_{\mathbf G}$-class $[x]_{\approx_{\mathbf G}}$ onto its corresponding $\approx_{\mathbf H}$-class $[y]_{\approx_{\mathbf H}}$ with $o(x)=o(y)$. Let us prove that $\vartheta$ is digraph isomorphism from $\vec{\mathcal G}(\mathbf G)$ to $\vec{\mathcal G}(\mathbf H)$. Suppose that $x_0\rightarrow_{\mathbf G}x_1$, and $\psi(x_0)=y_0$ and $\vartheta(x_1)=y_1$. Then $o(x_1)|o(x_0)$, and $x_0,x_1\in C_i^G$, for some $i\in n$. Since $y_0$ and $y_1$ are in corresponding $\approx_{\mathbf H}$-classes to $\approx_{\mathbf G}$-classes of $x_0$ and $x_1$, respectively, then $y_0,y_1\in C_i^H$. Since $y_0$ and $y_1$ are contained in a cyclic subgroup, and $o(y_1)|o(y_0)$, then $y_0\rightarrow_{\mathbf H}y_1$. This proves the theorem.
\end{dokaz}

\begin{pos}\label{sva tri su ista}
For all finite grous $\mathbf G$ and $\mathbf H$ the following conditions are equivalent :
\begin{enumerate}
\item Undirected power graphs of $\mathbf G$ and $\mathbf H$ are isomorphic;
\item Directed power graphs of $\mathbf G$ and $\mathbf H$ are isomorphic;
\item Enhanced power graphs of $\mathbf G$ and $\mathbf H$ are isomorphic.
\end{enumerate}
\end{pos}

The proof of the above corollary follows from Theorems \ref{neusmereni odredjuju usmerene} and  \ref{enhanced odredjuje usmereni}.

In \cite{power graph} Peter Cameron and Shamik Ghosh proved that finite abelian groups with isomorphic power graphs are isomorphic. Therefore, by Theorem 1 in \cite{power graph} and Corollary \ref{sva tri su ista} we get the following corollary.

\begin{pos}\label{na zanimljiviji nacin}
Let $\mathbf G$ and $\mathbf H$ be finite abelian groups such that $\mathcal G_e(\mathbf G)\cong\mathcal G_e(\mathbf H)$. Then $\mathbf G\cong\mathbf H$.%\footnote{Da li da napi\v semo na\v s dokaz, ne bi bio doga\v cak. Za nekog \v citaoca mo\v ze biti koristan jer Cameron nije kod sebe pokazao kako ta\v cno broj elemenata svakog reda odre\dj uje Abelovu grupu, a mi sad imamo to izvedeno na lep na\v cin.}
\end{pos}

%\begin{dokaz}
%Let $p$ be a prime with $p\big||G|$. Since $\mathbf G$ and $\mathbf H$ have the same numbers of elements of each order, then so do their $p$-Sylow subgroups $\mathbf G_p$ and $\mathbf H_p$. Now, each abelian $p$-group is isomorphic to $\prod_{i\in n}\mathbf C_{p^{j_i}}$, where $j_i>0$, and for each $j>0$ the group $\prod_{i\in n}\mathbf C_{p^{j_i}}$ has $\prod_{i\in n}p^{\min\{j,j_i\}}$ of elements of order at most $p^j$. This implies that abelian $p$-groups with the same numbers of elements of each order are isomorphic. It follows that the abelian groups $\mathbf G$ and $\mathbf H$ have isomorphic Sylow subgroups, and thus $\mathbf G\cong\mathbf H$.
%\end{dokaz}\\

In \cite{power graph} they also proved that, for any $k\in\mathbb N$, finite groups with isomorphic directed power graphs have the same numbers of elements of order $k$. Therefore, Proposition 4 in \cite{power graph} and Corollary \ref{sva tri su ista} give us the next corollary. %Further in this section, for any finite group and any $m\in\mathbb N$,  we are going directly from its enhanced power graph to determine how many elements of order $m$ does the group have.

\begin{pos}
Let $\mathbf G$ and $\mathbf H$ be finite groups with $\mathcal G_e(\mathbf G)\cong\mathcal G_e(\mathbf H)$. Then for every $m\in\mathbb N$ $\mathbf G$ and $\mathbf H$ have the same numbers of elements of order $m$.
\end{pos}

%The above follows from Proposition 4 in \cite{power graph} and Corollary \ref{sva tri su ista}.

Now, for any $m\in\mathbb N$ we will determine from the enhanced power graph how many elements of order $m$ does the group have.

Let $X$ be a finite set. For $A,B\in\mathcal P^+(X)$, we say $A\approx_{k}B$ if $k\big||A\cap B|$. We also define $\Phi_m:\mathcal P^+(X)\rightarrow \mathbb N$ in the following way:
\begin{align*}
\Phi_m(A)=\begin{cases}
\varphi(m)&,m\big||A|\\
0&,m\not{\big|}|A|,
\end{cases}
\end{align*}
where $\varphi$ is the Euler function. %The following lemma  follows from Lemma \ref{klike i ciklicne grupe}, and the fact that a cyclic group $n$ contains a unique subgroup of order $m$, for any $m|n$.

\begin{lem}\label{jeste ekvivalencija}
Let $\mathbf G$ be a finite group, and let $\mathcal Cl$ be the set of all maximal cliques of $\mathcal G_e(\mathbf G)$ of sizes divisible by $m$. Then $\approx_m$ is an equivalence relation on $\mathcal Cl$.
\end{lem}

\begin{dokaz}
Obviously $\approx_m$ is reflexive and symmetric relation. It remains to prove transitivity. Suppose that $C_1\approx_m C_2$ and $C_2\approx_m C_3$. By Lemma \ref{klike i ciklicne grupe} $\mathbf C_1$, $\mathbf C_2$, and $\mathbf C_3$ are cyclic subgroups of $\mathbf G$. Since $C_1\approx_m C_2$, there is a cyclic subgroup $\mathbf D_1\leq \mathbf C_1\cap \mathbf C_2$ of order $m$. Similarly there is a cyclic subgroup $\mathbf D_2\leq \mathbf C_2\cap\mathbf C_3$. Since the cyclic group $\mathbf C_2$ has a unique subgroup of order $m$, it follows that $\mathbf D_1=\mathbf D_2\leq C_1\cap C_3$. Therefore, $C_1\approx_m C_3$.
\end{dokaz}\\

For any group $\mathbf G$, and $X\subseteq G$, let us denote with $X^{(m)}$ the set of elements of order $m$ contained in $X$, and let us denote with $X^{(|m)}$ the set of elements in $X$ whose orders are divisors of $m$. Also, throughout this paper $\mathbb N_n$ shall denote the set $\{1,2,...,n\}$.

%The following is also a consequence of Corollary \ref{sva tri su ista} and Proposition 4 of \cite{power graph}. However, we still provide our own proof as it gives an easy way to calculate the numbers of elements of each order in a finite group from its enhanced power graph.

\begin{tvr}\label{prebrajanje}
Let $\mathbf G$ be a finite group. Let $\mathcal Cl=\{C_1,C_2,...,C_n\}$ be the set of all maximal cliques of $\mathcal G_e(\mathbf G)$, and let $\mathcal D=\{C\in\mathcal Cl\mid m\big||C|\}$. Then:
\begin{enumerate}
\item $|G^{(m)}|=\sum_{k=1}^n\big((-1)^{k+1}\sum_{I\subseteq \mathbb N_n,|I|=k}\Phi_m\big(\bigcap_{i\in I}C_i\big)\big)$;
\item $|G^{(m)}|=\varphi(m)\cdot\lvert\mathcal D/\approx_m\rvert$;
\end{enumerate}
\end{tvr}

\begin{dokaz}
{\it 1.} %For any finite group $\mathbf G$, let us denote with $G^{(m)}$ the set of all elements of $\mathbf G$ or order $m$. Also, for simpler notations we shall assume that $\mathbf G$ and $\mathbf H$ are finite groups such that $\Gamma=\mathcal G_e(\mathbf G)=\mathcal G_e(\mathbf H)$. 
By Lemma \ref{klike i ciklicne grupe}, a maximal clique of $\mathbf G_e(\mathbf G)$ is a maximal cyclic subgroup of $\mathbf G$. Therefore, $\mathbf C_1, \mathbf C_2, ..., \mathbf C_n$ are all maximal cyclic subgroups of $\mathbf G$. Since $G^{(m)}=\bigcup_{i=1}^n C_i^{(m)}$, we have the following:
\begin{align*}
|G^{(m)}|&=\Big|\bigcup_{i=1}^n C_i^{(m)}\Big|=\sum_{k=1}^n\Big((-1)^{k+1}\sum_{I\subseteq \mathbb N_n,|I|=k}\Big|\bigcap_{i\in I}C_i^{(m)}\Big|\Big)\\
&=\sum_{k=1}^n\Big((-1)^{k+1}\sum_{I\subseteq \mathbb N_n,|I|=k}\Big|\Big(\bigcap_{i\in I}C_i\Big)^{(m)}\Big|\Big).
\end{align*}
Now, for every $k\in\mathbb N$, a cyclic group of order $k$ has $\Phi(k)$ elements of order $m$, so we have:%, or none, depending of whether its order is divisible by $m$, or not. %(Here $\varphi$ is the Euler function.)
%
%By Lemma \ref{klike i ciklicne grupe}, it follows that if $C_0, C_1,..., C_{n-1}$ are the maximal cliques of $\Gamma$, and we have:
\begin{align*}
|G^{(m)}|=\sum_{k=1}^n\Big((-1)^{k+1}\sum_{I\subseteq \mathbb N_n,|I|=k}\Phi_m\Big(\bigcap_{i\in I}C_i\Big)\Big),
\end{align*}
which proves {\it 1}.

{\it 2.} By Lemma \ref{jeste ekvivalencija} the relation $\approx_m$ is an equivalence relation on $\mathcal D$. Let $\mathcal D_0$ be a $\approx_m$-class, and let $D_0\in \mathcal D_0$. Then $\mathbf D_0$ is a maximal cyclic subgroup of $\mathbf G$ with $m\big||D_0|$, so $\mathbf D_0$ has a unique cyclic subgroup $\mathbf M$ of order $m$. What is more, $\mathbf M$ is the unique cyclic subgroup of $\mathbf D$ for all $D\approx_m D_0$. This implies that all maximal cliques from $\mathcal D_0$ contain the same $\varphi(m)$ elements of order $m$, the generators of $\mathbf M$. Beside that, each element of order $m$ is contained in at least one maximal cyclic subgroup whose order is divisible by $m$, and consequently in at least one maximal clique from $\mathcal D$. Also, if there is an $x$ of order $m$ such that $x\in  E, F$ for $E,F\in \mathcal D$, then $\lb x\rb\leq\mathbf E,\mathbf F$, which implies $E\approx_m F$. This proves the equality.
\end{dokaz}\\

In the next proposition, given the enhanced power graph of an abelian group, we will give a description of this abelian group.

\begin{tvr}\label{prepoznavanje na osnovu prebrajanja}
Let $\mathbf G$ be a finite abelian group with $|G|=\prod_{i=1}^np_i^{k_i}$, where $p_i$ are all different prime numbers. For every $m\big||G|$ let $\mathcal Cl_m$ be the set of all maximal cliques of $\mathcal G_e(\mathbf G)$ with sizes divisible by $m$. Then:
\begin{align*}
\mathbf G\cong \prod_{i=1}^n\prod_{j=1}^{k_i}\big(\mathbf C_{p_i^{j}}\big)^{2L(p_i,j)-L(p_i,j-1)-L(p_i,j+1)}
\end{align*}
where $L(p,j)=\log_p\big(\sum_{k=0}^j \varphi(p^k)\cdot|\mathcal Cl_{p^k}/\approx_{p^k}|\big)$. 
\end{tvr}

\begin{dokaz}
Since $\mathbf G$ is an abelian group, $\mathbf G\cong \prod_{i=1}^n\mathbf G_{p_i}$, where, for each $i\in \mathbb N_n$, $\mathbf G_{p_i}$ is the unique $p_i$-Sylow subgroup of $\mathbf G$. Therefore, it is sufficient to prove the statement for any finite abelian $p$-group.

Let $\mathbf G$ be a finite abelian $p$-group. Then $\mathbf G\cong\prod_{i\in m}\mathbf C_{p^{l_i}}$, where $l_i>0$ for all $i\in n$. We shall call the cyclic groups $\mathbf C_{p^{l_i}}$ factors of $\mathbf G$. It is easy to see that \begin{align*}
\big|G^{(|p^j)}\big|=\prod_{i\in m}p^{\min\{j,l_i\}}=p^{\sum_{i\in m}\min\{j,l_i\}},
\end{align*}
for any $j\geq 0$. Therefore, for any $j\geq 1$, there are $\sum_{i\in m}\min\{j,l_i\}-\sum_{i\in m}\min\{j-1,l_i\}=\log_p (|G^{(|p^j)}|)-\log_p (|G^{(|p^{j-1})}|)$ factors of $\mathbf G$ whose orders are at least $p^j$, which is by Proposition \ref{prebrajanje} equal to $L(p,j)-L(p,j-1)$. Now, the number of  factors of $\mathbf G$ whose orders are exactly $p^j$ is $L(p,j)-L(p,j-1)-(L(p,j+1)-L(p,j))=2L(p_i,j)-L(p_i,j-1)-L(p_i,j+1)$. This proves the proposition.
\end{dokaz}\\

Note that there are non-isomorphic groups with same numbers of elements of each order. Therefore, given $\mathcal G_e(\mathbf G)$, if we do not know whether $\mathbf G$ is abelian, we may not be able to recognize the group $\mathbf G$. Pair of groups of order $27$ given by Cameron and Ghosh are one such example: $\mathbf C_3\times\mathbf C_3\times\mathbf C_3$ and the group with presentation $\mathbf G=\langle x,y\mid x^3=y^3=[x,y]^3=1  \rangle$.

\section{Automorphism groups of enhanced power graphs}
\label{grupe automorfizama}

As we have seen, power graph of a finite group determines the directed power graph of the group. It is not hard to see that an automorphism of power graph of a group may not be an automorphism of directed power graph of a group. However, as we shall prove, an automorphism of power graph of a finite group is an automorphism of enhanced power graph of the group.

\begin{teo}\label{power indukuje enhanced}
Let $\mathbf G$ and $\mathbf H$ be finite groups. Then any isomorphism from $\mathcal G(\mathbf G)$ to $\mathcal G(\mathbf H)$ is an isomorphism from $\mathcal G_e(\mathbf G)$ to $\mathcal G_e(\mathbf H)$.
\end{teo}

\begin{dokaz}
Here we shall refer to some results of Cameron from \cite{power graph 2}.

If $\mathbf G$ is a generalized quaternion group, then its enhanced power graph and power graph are equal since it is a $p$-group, implying $\Aut(\mathcal  G(\mathbf G))=\Aut(\mathcal  G_e(\mathbf G))$. Also, if $\mathbf G$ is a cyclic group, then its enhanced power graph is a complete graph, and then we obviously have $\Aut(\mathcal  G(\mathbf G))\subseteq\Aut(\mathcal  G_e(\mathbf G))$.

So for the rest of the proof we shall assume that $\mathbf G$ is neither generalized quaternion group, nor cyclic group. By Proposition 4 of \cite{power graph 2} the identity element of $\mathbf G$ is the only one being adjacent to all other elements of $\mathcal G(\mathbf G)$. %In this case Cameron has introduced the following notions:
%\begin{itemize}
%\item Closed neighborhood of $x$, denoted with $\overline N(x)=\{y\in G\mid x\sim y\text{ or }x=y\}$;
%\item Relation $\equiv$ defined with $x\equiv y$ if $\overline N(x)=\overline N(y)$;
%\item Relation $\approx$ defined with $x\approx y$ if $\lb x\rb=\lb y\rb$.
%\end{itemize}

Both relations $\equiv$ and $\approx$ are easily seen to be equivalence relations, with $\equiv$-classes being recognizable in the graph. Moreover, in Proposition 5 of \cite{power graph 2} Cameron showed what is the relation between the two equivalences. $\approx$ is a refinement of $\equiv$, and each non-identity $\equiv$-class $C$ is of one of the following forms:
\begin{enumerate}
\item[$(1.)$] $C$ is a $\approx$-class;
\item[$(2.)$]  $C=\{x\in\lb y\rb\mid o(x)\geq p^s\}$, where $o(y)=p^r$, $1\leq s\leq r-1$, and $p$ is a prime number.
\end{enumerate}
Further in his proof he also shows that we can recognize the $\equiv$-classes of type $(2.)$ in the power graph, and recognize the prime number related to that $\equiv$-class.

%Next he shows that we can reconstruct the directed power graph up to isomorphism in the following way:
%\begin{itemize}
%$\item For any pair of elements which do not belong to a same $\equiv$-class of type $(2.)$, we know exactly how they relate in the directed power graph.
%\item We know that every $\equiv$-class of type $(2.)$ is consisted of $r-s+1$ $\approx$-classes of elements of orders $p^s,p^{s+1},...,p^r$ in the group.
%\end{itemize}
%But the only arrows of $\vec{\mathcal G}(\mathbf G)$ that cannot exactly be reconstructed from $\mathcal G(\mathbf G)$, but only up to isomorphism, are those connecting elements belonging to a same $\equiv$-class of type $(2.)$.
%
%We are going to show that using the above said we can recognize in $\mathcal G(\mathbf G)$ maximal cliques of $\mathcal G_e(\mathbf G)$, which are  by Lemma \ref{klike i ciklicne grupe} maximal cyclic subgroups of $\mathbf G$. Indeed, $x$ and $y$ belong to a maximal cyclic subgroup of $\mathbf G$ iff they belong to any cyclic subgroup of $\mathbf G$, i.e. there exists an element $z\in G$ so that $z\rightarrow x$ and $z\rightarrow y$.\\

For any $x\in G$ the size of $[x]_\approx$ is $\varphi(o(x))$, where $\varphi$ is the Euler function. It is also known that $\varphi(\prod_{j=1}^l p_j^{i_j})=\prod_{j=1}^l (p_j^{i_j-1}(p_j-1))$ for any prime numbers $p_1,p_2,...,p_l$ and any $i_1,i_2,...,i_l>0$. Therefore, $n|m$ implies $\varphi(n)|\varphi(m)$, while equality only holds when $n=m$, or $n$ is odd and $m=2n$. This implies that, for any $x$ and $y$, $x\rightarrow y$ if and only if:
\begin{align}
&x\sim y,%\ \big|[y]_\approx\big|\Big|\big|[x]_\approx\big|,
\text{ and } \big|[y]_\approx\big|<\big|[x]_\approx\big|,\text{ or}\label{prvi slucaj}\\
\begin{split}
&x\sim y,\  \big|[y]_\approx\big|=\big|[x]_\approx\big|, \text{ and }x\sim z\text{ for some } z\text{ such that}\\
&\ \ \ \ \ \ \ \ \ \ \ \ \ \ \ \ \ \ \ \ \ \ \ \ \ \ \ \ \ \ \    [z]_\approx=\{z\},\text{ and } \overline N(z)\neq G,\text{ or}\label{drugi slucaj}
\end{split}\\
&x\approx y\label{treci slucaj}.
\end{align}
Note that $z$ in the condition (\ref{drugi slucaj}) is an involution, so this one covers the case when $\big|[x]_\approx\big|=2\big|[y]_\approx\big|$ and $\big|[y]_\approx\big|$ is an odd number. In this case $[z]_\approx=[z]_\equiv$.

It is easily seen that if $x_0\rightarrow y_0$, then $x\rightarrow y$ for any $x\approx x_0$ and $y\approx y_0$. This is true even for $\equiv$-classes, unless $x_0$ and $y_0$ are in the same class of type $(2.)$. Suppose that $[y_0]_\equiv$ is of type $(2.)$, and $x_0\rightarrow y_0$ and $x_0\not\equiv y_0\equiv y$. Then $x_0\sim y$, i.e. $x_0\rightarrow y$ or $y\rightarrow x_0$. But $y\rightarrow x_0$ and $x_0\rightarrow y_0$ would imply that $x_0\in[y_0]_\equiv$. We have shown that $x_0\rightarrow y_0$ and $x_0\not\equiv y_0\equiv y$ implies $x_0\rightarrow y$. %It it has to be $x_0\rightarrow y$.
%
%Indeed, if $x_0\rightarrow y_0$ holds for some $y_0$ in a $\equiv$-class of type $(2.)$, and $x_0\not\in[y_0]_\equiv$, then $p^r|o(x_0)$, and $p^r<o(x_0)$. Otherwise $x_0\in[y_0]_\equiv$ would hold. Now for any $y\equiv y_0$ we have $x_0\sim y$, and because $p^r<o(x_0)$ we have $x_0\rightarrow y$.
%
In similar way can be proved that $x_0\rightarrow y_0$ and $x\equiv x_0\not\equiv y_0$ implies $x\rightarrow y_0$.
%
%Also, if we have $x_0\rightarrow y_0$ for some $x_0$ in a $\equiv$-class of type $(2.)$, and $y_0\not\in[x_0]_\equiv$, then $o(y_0)=p^t$ for some $t<s$, and this implies $x\rightarrow y_0$ for all $x\equiv x_0$.
%
From this we can easily conclude $x_0\rightarrow y_0$ implies $x\rightarrow y$ for all $x\equiv x_0$ and $y\equiv y_0$, whenever $x_0$ and $y_0$ are not in a same $\equiv$-class of type $(2.)$.

%Now, for power graph of $\mathcal G(\mathbf G)$ of finite group $\mathbf G$ we introduce graph $\mathcal E(\mathbf G)$ with $V(\mathcal E(\mathbf G))=G$, and $x\sim_{\mathcal E(\mathbf G)} y$ iff there is $z\in G$ so that:
%\begin{itemize}
%\item $z$ and $x$ satisfy either (\ref{prvi slucaj}), or (\ref{drugi slucaj}), or (\ref{treci slucaj}), and
%\item $z$ and $y$ satisfy either (\ref{prvi slucaj}), or (\ref{drugi slucaj}), or (\ref{treci slucaj}).
%\end{itemize}

Thanks to the above said for any $x,y\in G$, whenever $x\not\equiv y$, we can decide from $\mathcal G(\mathbf G)$ whether $x\rightarrow y$. Now we prove that for any $x,y\in G$ we can decide whether $x\sime y$. Indeed, if $x\sim y$, then $x\sime y$. On the other hand, if $x\not\sim y$, then $x\sime y$  if and only if  there exists $z\in G$ for which $z\rightarrow x$, $z\rightarrow y$, and $x\not\equiv y\not\equiv z\not\equiv x$. Given that, for any finite group $\mathbf G$, it is decidable from $\mathcal G(\mathbf G)$ whether $x\sime y$.

Now, let $\mathbf G$ and $\mathbf H$ be finite groups, and let $\varphi:G\rightarrow H$ be a graph isomorphism from $\mathcal G(\mathbf G)$ to $\mathcal G(\mathbf H)$. Since whether $x\sime_{\mathbf G} y$ is decidable from $\mathcal G(\mathbf G)$, and whether $\varphi(x)\sime_{\mathbf H}\varphi(y)$ is decidable from $\mathcal G(\mathbf H)$, it follows that $\varphi$ is also a graph isomorphism from $\mathcal G_e(\mathbf G)$ to $\mathcal G_e(\mathbf H)$.
%
%Thanks to the above said now we can transform all power graph edges bet-ween different $\equiv$-classes into directed power graph arrows based on conditions (\ref{prvi slucaj}) and (\ref{drugi slucaj}), while we would take the $\equiv$-classes of type $(2.)$ as if they had $\varphi(p^r)$ elements (with $p$ being the specific prime number for each such $\equiv$-class). Now for any $x$ and $y$ in a same $\equiv$ class we will take $x\sime y$. Next, for any $x$ and $y$, for which there is some $z$ so that $z\rightarrow x$ and $z\rightarrow y$, we also take $x\sime y$. And finally, for any $x$ and $y$ so that $x\rightarrow y$ we take $x\sime y$. This way we have reconstructed the entire $\mathcal G_e(\mathbf G)$.
%
%Indeed, if $x$ and $y$ are adjacent in $\mathcal G_e(\mathbf G)$, then there is an element $z$ such that $z\rightarrow x$ and $z\rightarrow y$, or $x\rightarrow y$, or $y\rightarrow x$ in  $\vec{\mathcal G}(\mathbf G)$. If $x\approx y$, then $x\equiv y$, and we do get $x\sime y$ by the above procedure. If $x\not\approx y$, and $x\rightarrow y$ or $y\rightarrow x$, then we also get $x\sime y$. And if $z\rightarrow x$ and $z\rightarrow y$ for $x\not\approx y\not\approx z\not\approx x$, then we again get $x\sime y$ by the above procedure.
%
%We have proven that from $\mathcal G(\mathbf G)$ we can reconstruct exactly $\mathcal G_e(\mathbf G)$, and not just up to isomorphism. This implies that $\Aut(\mathcal G(\mathbf G))\subseteq\Aut(\mathcal G_e(\mathbf G))$.
\end{dokaz}

\begin{pos}\label{odnos automorfizama}
Let $\mathbf G$ be a finite group. Then
$$\Aut(\mathbf G)\subseteq\Aut(\vec{\mathcal  G}(\mathbf G))\subseteq \Aut(\mathcal  G(\mathbf G))\subseteq\Aut(\mathcal  G_e(\mathbf G)).$$
\end{pos}

\begin{dokaz}
One can easily see that $\Aut(\mathbf G)\subseteq\Aut(\vec{\mathcal  G}(\mathbf G))\subseteq \Aut(\mathcal  G(\mathbf G))$, and the inclusion $ \Aut(\mathcal  G(\mathbf G))\subseteq\Aut(\mathcal  G_e(\mathbf G))$ is a direct consequence of Theorem \ref{power indukuje enhanced}.
\end{dokaz}\ \\

Note that in the above theorem those inclusions for some groups can also be strict. For example $\Aut(\mathbf C_5)\subsetneq\Aut(\vec{\mathcal  G}(\mathbf C_5))\subsetneq \Aut(\mathcal  G(\mathbf C_5))$, and $\Aut(\mathcal  G(\mathbf C_6))\subsetneq\Aut(\mathcal  G_e(\mathbf C_6))$, so no equality holds in general.

\begin{pos}
The Klein group is the only non-trivial finite group $\mathbf G$ such that $\Aut(\mathbf G)=\Aut\big(\mathcal G_e(\mathbf G)\big)$.
\end{pos}

\begin{dokaz}
Suppose that $\mathbf G$ is a finite non-trivial group with $\Aut(\mathbf G)=\Aut(\mathcal G_e(\mathbf G))$. By Theorem \ref{odnos automorfizama} it follows $\Aut(\mathbf G)=\Aut(\mathcal G(\mathbf G))$, which implies $\mathbf G\cong \mathbf C_2\times\mathbf C_2$. On the other hand, it is easily checked that $\Aut (\mathbf C_2\times\mathbf C_2)=\Aut(\mathcal G(\mathbf C_2\times\mathbf C_2))$.
\end{dokaz}

\begin{lem}\label{automorfizmi u jakom proizvodu}
If $\Gamma$ and $\Delta$ are any two graphs, then $\Aut(\Gamma)$ can be embedded into $\Aut(\Gamma\boxtimes\Delta)$.
\end{lem}

\begin{dokaz}
Let us define the mapping $\Phi:\varphi\mapsto\hat\varphi$ with $\hat\varphi\big((x,y)\big)=(\varphi(x),y)$ for all $\varphi\in\Aut(\Gamma)$. Let $\varphi_0\in\Aut(\Gamma)$, and let $\hat\varphi_0=\Phi(\varphi_0)$. Obviously $\hat\varphi_0$ is a bijection. Also, for all $x_1,x_2\in V(\Gamma)$ and $y_1,y_2\in V(\Delta)$ we have $(x_1,y_1)\sim(x_2,y_2)$  if and only if  $(\varphi_0(x_1),y_1)\sim(\varphi_0(x_2),y_2)$  if and only if  $\hat\varphi_0((x_1,y_1))\sim\hat\varphi_0((x_2,y_2))$, so $\hat\varphi_0\in\Aut(\Gamma\boxtimes\Delta)$. Besides, the implication $\varphi\neq\psi\Rightarrow\Phi(\varphi)\neq\Phi(\psi)$ obviously holds, and it is easily seen that $\Phi:\Aut(\Gamma)\rightarrow\Aut(\Gamma\boxtimes\Delta)$ is a homomorphism. Therefore, $\Phi$ is an injective homomorphism, which proves the statement of the lemma.
\end{dokaz}

\begin{teo}
$\mathbf C_2$ is the only non-trivial finite group whose enhanced power graph has abelian automorphism group. 
\end{teo}

\begin{dokaz}
Let $\mathbf G$ be a non-trivial finite group, such that $\Aut(\mathcal G_e(\mathbf G))$ is abelian. If $\mathbf G$ is not nilpotent, then $\mathbf G/\Z(\mathbf G)$ is not abelian, and neither is $\Aut(\mathcal G_e(\mathbf G))$ since $\mathbf G/\Z(\mathbf G)\cong\Inn(\mathbf G)\leq \Aut(\mathcal G_e(\mathbf G))$.

Since $\mathbf G$ is nilpotent, by Theorem 5.39 in \cite{rotman} its Sylow subgroups are unique for all prime divisors of $|\mathbf G|$, and $\mathbf G$ is their direct product. By Lemma \ref{proizvodi} %it follows that $\mathcal G_e(\mathbf G)=\mathcal G_e(\mathbf P_1)\boxtimes\mathcal G_e(\mathbf P_2)\boxtimes ...\boxtimes\mathcal G_e(\mathbf P_k)$. Now, if an $\Aut(\mathcal G_e(\mathbf P_i))$ is non-abelian, then so is $\Aut(\mathcal G_e(\mathbf G))$. Indeed, if $\varphi,\psi\in\Aut(\mathcal G_e(\mathbf P_i))$ don't commute, then we define $\hat\varphi$ and $\hat\psi$ with 
%\begin{align*}
%&\hat\varphi\big((x_1,...,x_k)\big)=(x_1,...,x_{i-1},\varphi(x_i),x_{i+1},...,x_k),\text{ and}\\
%&\hat\psi\big((x_1,...,x_k)\big)=(x_1,...,x_{i-1},\psi(x_i),x_{i+1},...,x_k).
%\end{align*}
%It is easy to see that $\hat\varphi,\hat\psi\in\Aut(\mathcal G_e(\mathbf G))$, and that $\hat\varphi$ and $\hat\psi$ don't commute. So 
 and Lemma \ref{automorfizmi u jakom proizvodu} it is sufficient to prove that the only non-trivial $p$-group with abelian $\mathcal G_e(\mathbf G)$ is $\mathbf C_2$.

%If both $\mathbf P_2$ and $\mathbf P_3$ are cyclic, then $\mathbf G$ is also cyclic of size at least 6, and $\mathcal G_e(\mathbf G)$ is non-abelian. If at least one of $\mathbf P_2

%We know that for the 3-group $\mathbf P_3$ $\Aut(\mathcal G_e(\mathbf P_3))$ has a subgroup $\mathbf F$ isomorphic to $\mathbf S_3$. Now, each elements $\varphi$ of $\mathbf F$ can be extended from $\mathbf P_3$ to $\mathbf G$ with $\hat\varphi:(x,y)\mapsto(x,\varphi(y))$, for all $x\in P_2$ and $y\in P_3$. It is easily seen that $\hat\varphi$ is a bijection on $G$. Furthermore, $\hat\varphi$ defined like this is an automorphism on $\mathcal G_e(\mathbf G)$. By Lemma \ref{proizvodi} it is enough to prove that $\hat\varphi$ is an automorphism on $\Gamma=\mathcal G_e(\mathbf P_2)\boxtimes\mathcal G_e(\mathbf P_3)$, and we do have $(x_1,y_1)\sim_\Gamma(x_2,y_2)$, iff $x_1\sime_{\mathbf P_2}x_2$ and $y_1\sime_{\mathbf P_3}y_2$, iff $\varphi(x_1)\sime_{\mathbf P_2}\varphi(x_2)$ and $y_1\sime_{\mathbf P_3}y_2$, iff $\hat\varphi(x_1,y_1)\sim_\Gamma\hat\varphi(x_2,y_2)$. What is more, if $\varphi_1\neq\varphi_2$ implies $\hat\varphi_1\neq\hat\varphi_2$, so $\hat{\mathbf F}$ with universe $\hat F=\{\hat\varphi\mid\varphi\in F\}$ is a subgroup of $\Aut(\mathcal G_e(\mathbf G))$ such that $\mathbf F\cong\hat{\mathbf F}$. Now, since $\Aut(\mathcal G_e(\mathbf G))$ has non-abelian subgroup, $\Aut(\mathcal G_e(\mathbf G))$ is not abelian either. This proves the statement of the theorem.

If a $p$-group $\mathbf G$ has an element of order at least 5, then there are $x,y,z\in G$ so that $\lb x\rb=\lb y\rb=\lb z\rb$. Then permutations $(x,y)$ and $(y,z)$ on $G$ are automorphisms of $\mathcal G_e(\mathbf G)$ which do not commute.
%is a $p$-group for a $p\geq 5$, then there is an element $x\in G$ so that $o(x)\geq 5$ is a prime number. Then we have $\mathbf S_3\cong\lb(x,x^2),(x^2,x^3)\rb\leq\Aut(\mathcal G_e(\mathbf G))$, thus $\Aut(\mathcal G_e(\mathbf G))$ is non-abelian.
Let $\mathbf G$ be a $p$-group of exponent 3. %If it contains an element of order $9$, then $\mathbf S_3\cong\lb (x,x^2),(x^2,x^4)\rb\leq \Aut(\mathcal G_e(\mathbf G))$. So Let us assume that the exponent of $\mathbf G$ is 3.
If $\mathbf G\cong\mathbf C_3$, then $\Aut(\mathcal G_e(\mathbf G))\cong \mathbf S_3$. If it is non-cyclic, then $|\mathbf G|\geq 9$ and $\mathbf G$ has exponent $3$, and there are elements $x$, $y$, and $z$ that generate three different cyclic subgroups. Then $(x,y)(x^2,y^2),(y,z)(y^2,z^2)\in \Aut(\mathcal G_e(\mathbf G))$, and they do not commute. Now let $\mathbf G$ be a 2-group of exponent 4. 
%If its exponent is at least 8, then there is an element $x$ of order 8, and $S_3\cong\lb(x,x^3),(x^3,x^5)\rb\leq\Aut(\mathcal G_e(\mathbf G))$. Let us assume now that the exponent of $\mathbf G$ is 4. 
If there is an element $\mathbf G$ of order 2 that belongs to exactly one cyclic subgroup $\lb x\rb$ of order 4, then $(x,x^2),(x^2,x^3)\in\Aut(\mathcal G_e(\mathbf G))$, and $(x,x^2)$ and $(x^2,x^3)$ do not commute. If there is no such element of order 2, then there are $x,y\in G$ of order 4 so that $|\lb x\rb\cap\lb y\rb|=2$. Then $(x,x^3)$ and $(x,y)(x^3,y^3)$ are automorphisms of $\mathcal G_e(\mathbf G)$ that do not commute. Finally, it is easy to notice that $\mathbf C_2$ is the only group of exponent 2 having abelian automorphism group of its enhanced power graph. This proves the theorem.
%Let $\mathbf G$ be a group with exponent $2$. If $\mathbf G\cong \mathbf C_2$, then $\Aut(\mathcal G_e(\mathbf G))\cong\mathbf C_2$ is an abelian group. If however $\mathbf G$ is not cyclic, then it has at least 3 elements $x,y,z\in G$ of order 2, and $\mathbf S_3\cong \lb(x,y),(y,z)\rb\leq\Aut(\mathcal G_e(\mathbf G))$.
%
%Let $|\mathbf G|=2^{n_1}3^{n_2}$ for $n_1,n_2>0$. 
\end{dokaz}

\begin{teo}
$\mathbf C_2$, $\mathbf C_3$, and the Klein group are the only non-trivial finite groups whose automorphism groups of their enhanced power graphs have square free orders.
\end{teo}

\begin{dokaz}
Let $\mathbf G$ be a group with $|\Aut(\mathcal G_e(\mathbf G))|$ being square free. If $\mathbf G$ is cyclic, then $\mathbf G\cong\mathbf C_2$ or $\mathbf G\cong\mathbf C_3$, so assume that $\mathbf G$ is non-cyclic. Then for all $x\in G$ we have $o(x)<5$, because $o(x)\geq 5$ implies $\mathbf C_2\times\mathbf C_2\cong\lb(x,x^{-1}),(x^2,x^{-2})\rb\leq \Aut(\mathcal G_e(\mathbf G))$. Also, $\mathbf G$ has at most two elements of order larger than 2. Otherwise there would be $x$ and $y$ with $o(x),o(y)>2$ and $\langle x\rangle\neq\langle y\rangle$, and then $\mathbf C_2\times\mathbf C_2\cong\langle(x,x^{-1}),(y,y^{-1})\rangle\leq\Aut(\mathcal G_e(\mathbf G))$. If the non-cyclic group $\mathbf G$ had exactly one cyclic subgroup of order 3 or 4, then it would have at least three more involutions, and $\mathbf C_2\times\mathbf C_2$ could again be embedded into $\Aut(\mathcal G_e(\mathbf G))$. So the exponent of $\mathbf G$ is 2. Because $\mathbf G$ does not have more than 3 involutions, it follows $\mathbf G\cong \mathbf C_2\times\mathbf C_2$.
\end{dokaz}

\begin{teo}
$\mathbf C_2$ and $\mathbf C_4\times\mathbf C_2$ are the only non-trivial finite groups whose automorphism groups of their enhanced power graphs have prime power orders.
\end{teo}

\begin{dokaz}
Let $\mathbf G$ be a non-trivial finite group for which $|\Aut(\mathcal G_e(\mathbf G))|$ is a prime power. Since $\Aut(\mathcal G_e(\mathbf G))$ contains a transposition on $G$, $\Aut(\mathcal G_e(\mathbf G))$ is a 2-group.

If $\mathbf G$ is not nilpotent, then $\mathbf G/Z(\mathbf G)$ is not a $p$-group, and because $\mathbf G/Z(\mathbf G)\cong \Inn(\mathbf G)\leq \Aut(\mathcal G_e(\mathbf G))$, $\Aut(\mathcal G_e(\mathbf G))$ is not a $p$-group either. It follows that $\mathbf G$ is a nilpotent group. But no Sylow subgroup of $\mathbf G$ contains an element $x$ with $o(x)\geq 5$. Otherwise there would be $x,y,z\in G$ so that $\lb x\rb=\lb y\rb=\lb z\rb$ and $x\neq y\neq z\neq x$, and then the permutation $(x,y,z)$ would  be an automorphism of $\mathcal G_e(\mathbf G)$. It follows by  Theorem 5.39 in \cite{rotman} that $\mathbf G=\mathbf P_2\times\mathbf P_3$, where $\mathbf P_2$ and $\mathbf P_3$ are the 2-Sylow subgroup of $\mathbf G$ whose exponent is less than 8, and the 3-Sylow subgroup of $\mathbf G$ whose exponent is less than 9, respectively. Next, by Lemma \ref{proizvodi} we get $\mathcal G_e(\mathbf G)=\mathcal G_e(\mathbf P_2)\boxtimes\mathcal G_e(\mathbf P_3)$.

To show that $\mathbf P_3$ is trivial, it is sufficient to prove that otherwise $\Aut(\mathcal G_3(\mathbf P_3))$ would contain an element of order 3, because by Lemma \ref{automorfizmi u jakom proizvodu} that implies that $\Aut(\mathcal G_e(\mathbf G))$ also contains an element of order 3.  If $\mathbf P_3\cong \mathbf C_3$, then $\Aut(\mathcal G_e(\mathbf P_3))\cong\mathbf S_3$. If the $\mathbf P_3$ of exponent 3 is non-cyclic, then it has three elements $x$, $y$, and $z$ generating three different cyclic subgroups, and then $(x,y,z)(x^2,y^2,z^2)\in\Aut(\mathcal G_e(\mathbf G))$. Therefore, $\mathbf G$ is a 2-group.

Notice that $\mathbf G$ has at most 2 maximal cyclic subgroups of order two, and each non-maximal cyclic subgroup of order 2 is contained in exactly 2 maximal cyclic subgroups of order 4. What is more, group $\mathbf G$ contains at most two non-maximal cyclic subgroups of order 2. Otherwise $\mathcal G_e(\mathbf G)$ would have an automorphism of order 3. It follows that $|G|\leq 8$. Now it is easy to see that $\mathbf G$ is isomorphic to $\mathbf C_2$ or $\mathbf C_4\times\mathbf C_2$, and we have $\Aut(\mathcal G_e(\mathbf C_2))\cong \mathbf C_2$ and $\Aut(\mathcal G_e(\mathbf C_4\times\mathbf C_2))\cong\mathbf C_2\times\mathbf C_2\times\mathbf C_2\times\mathbf C_2$.
\end{dokaz}

%\section{Nesvrstani}

\section{The structure of enhanced power graphs of finite abelian groups}\label{prepoznavanje grupe}
\label{p semistabla}

%Our motivation in this section is to find a way to decide whether a given graph is enhanced power graph of some finite abelian group.
In this section we will describe enhanced power graphs of finite abelian groups. For start we will describe enhanced power graphs of finite abelian $p$-groups. To that end we introduce notions of rooted $p$-trees and $p$-semitrees. We will prove that $p$-semitrees are exactly enhanced power graphs of finite abelian $p$-groups.

A {rooted tree} is a tree whose one vertex is designated as the root. For two vertices $x$ and $y$ of a rooted tree we write $x<y$ if the unique path from $y$ to the root passes through $x$. If $x$ and $y$ are adjacent and $x<y$, then $y$ is a child of $x$, and $x$ is the parent of $y$. The height %\footnote{Da li da koristim depth, ili height? U na\v sem slu\v caju je height logi\v cnije, ali se taj pojam izgleda uobi\v cajeno zove dubina, a visina, tj. height, je rastojanje od najudaljenijeg lista iznad njega. Prema Wikipediji. Sad sam ipak stavio height.}
  of vertex $x$ is the distance between $x$ and the root.

For tuples $\overline x=(x_1,x_2,...,x_n)$ and $\overline y=(y_1,y_2,...,y_n)$ of non-negative integers we say $\overline x\leq \overline y$ if $\overline x(i)\leq \overline y(i)$ for all $i\leq n$. Let $\overline a=(a_1,a_2,...,a_n)$ be a tuple of positive integers. Let $\overline b=(b_1,b_2,...,b_n)$ and $\overline c=(c_1,c_2,...,c_n)$ be tuples of integers  with $\overline b,\overline c\in\overline a^{\downarrow}=\{\overline x\mid \overline 0\leq \overline x\leq \overline a\}$. Then $\overline a$-width of the tuple $\overline b$, or simply width of $\overline b$, is the number $|\{i\in \mathbb N_n\mid b_i\neq a_i\}|$. We shall denote it with $w_{\overline a}(\overline b)$, or simply $w(\overline b)$. $\overline a$-height of tuple $\overline b$, or simply height of $\overline b$, is the number $\max\{a_i-b_i\mid i\in \mathbb N_n\}$. We shall denote it with $h_{\overline a}(\overline b)$, or simply $h(\overline b)$. We say that $\overline c$ is an $\overline a$-successor of $\overline b$, or simply a successor of $\overline b$, if $\overline c\neq\overline a$, and:
\begin{itemize}
\item $c_i=b_i-1$ if $b_i\neq a_i$, and
\item $c_i\in\{a_i,a_i-1\}$  if $b_i=a_i$.
\end{itemize}
If $\overline c$ is a successor of $\overline b$, then $\overline b$ is a predecessor of $\overline c$.

 The {\bf rooted $p$-tree} associated to a tuple $\overline a=(a_1,a_2,...,a_n)$ of positive integers, denoted with $T_p(\overline a)$, is the rooted tree $T$ whose vertices can be labeled with tuples from $\overline a^{\downarrow}$ in the following way (here we allow more vertices to be labeled with a same tuple):
\begin{itemize}
\item The root of $T$ is labeled with $\overline a$;
\item The root has exactly $\frac{p^n-1}{p-1}$ children which are all labeled with successors of $\overline a$ in the following way: For each successor $\overline b$ of $\overline a$ there are exactly $(p-1)^{w(\overline b)-1}$ children of the root labeled with $\overline b$;
\item If a vertex $x$ of $T$ is labeled with $\overline c\neq\overline a$, with $c_i\neq 0$ for all $i\in n$, then it has exactly $p^{n-1}$ children which are all labeled with successors of $\overline c$ in the following way: For each successor $\overline d$ of $\overline c$ there are exactly $p^{w(\overline c)-1}(p-1)^{w(\overline d)-w(\overline c)}$ children of $x$ labeled with $\overline d$;
\item If a vertex $x$ of $T$ is labeled with $\overline c\neq\overline a$, with $c_i= 0$ for some $i\leq n$, then $x$ has no children.
\end{itemize}
It is easily seen that any two rooted trees satisfying the above conditions are isomorphic, so $T_p(\overline a)$ is well defined.

The {\bf $p$-semitree} associated to the tuple $\overline a$, denoted by $S_p(\overline a)$, is the graph constructed from $T_p(\overline a)$ by adding an edge between every two $x$ and $y$ for which $x<y$, and by replacing each vertex of height $k>0$ by $K_{p^{k-1}(p-1)}=K_{\varphi(p^k)}$.

\begin{tvr}\label{karakterizacija}
A graph is the enhanced power graph of a finite abelian $p$-group  if and only if  it is a $p$-semitree.
\end{tvr}

\begin{dokaz}
The theorem obviously holds if the dimension of the tuple $\overline a$ is $1$, because, for any $m\in\mathbb N$, $\overline a=(m)$  if and only if  $S_p(\overline a)\cong K_{p^m}$, which is the enhanced power graph of $\mathbf C_{p^k}$. So we assume its dimension is at least $2$.

Let $\overline a=(a_1,a_2,...,a_n)$ be a tuple of positive integers, and let $\mathbf G=\prod_{i=1}^n \mathbf C_{p^{a_i}}$. Since every finite abelian abelian group is isomorphic to a direct product of cyclic groups, it is sufficient to prove that $S_p(\overline a)\cong \mathcal G_e(\mathbf G)$. It is easily seen that $S_p(\overline a)\cong\mathcal G_e(\mathbf G)$  if and only if  $T_p(\overline a)$ is, as a poset, isomorphic to the poset $P$ of cyclic subgroups of $\mathbf G$. To show that there is an isomorphism $\psi:T_p(\overline a)\rightarrow P$ which maps each vertex labeled with $\overline c=(c_1,c_2,...,c_n)$ into a cyclic subgroup whose $i$-th projection has order $p^{a_i-c_i}$ for each $i\leq n$, it is sufficient to prove the following:
%
%Let us regard each vertex of labeled $T_p(\overline a)$ labeled with $\overline c$ as a cyclic subgroup of $\mathbf G=\prod_{i\in n} \mathbf C_{p^{a_i}}$, whose each $i$-th projection is cyclic group of order $p^{a_i-c_i}$, naturally, with $x<y$ in the rooted tree if ; then it is enough to prove the following:
\begin{enumerate}
\item[$(1)$] $\mathbf G$ has exactly $\frac{p^n-1}{p-1}$ cyclic subgroups of order $p$;
\item[$(2)$] For $M\subseteq \mathbb N_n$, there are exactly $(p-1)^{|M|-1}$ cyclic subgroups of $\mathbf G$ of order $p$ whose $i$-th projection is non-trivial  if and only if  $i\in M$;
\item[$(3)$] If $\mathbf H$ is a cyclic subgroup of $\mathbf G$, and $\mathbf H_i$ its $i$-th projection for all $i\leq n$, then $\mathbf H$ is a maximal cyclic subgroup of $\mathbf G$  if and only if  $|H_i|=p^{a_i}$ for some $i\leq n$.
\item[$(4)$] If $\mathbf H$ is a non-maximal cyclic subgroup of $\mathbf G$ of order $p^h>1$ whose $i$-th projection is $\mathbf H_i$ of order $p^{h_i}$ for all $i\leq n$, then $\mathbf H$ is subgroup of exactly $p^{n-1}$ cyclic subgroups of $\mathbf G$ of order $p^{h+1}$. What is more, if $L=\{i\in \mathbb N_n\mid h_i\neq 0\}$ and $L\subseteq M\subseteq \mathbb N_n$, then $\mathbf H$ is subgroup of exactly $p^{|L|-1}(p-1)^{|M|-|L|}$ cyclic subgroups of order $p^{h+1}$ whose $i$-th projection is non-trivial  if and only if  $i\in M$.
\end{enumerate}

$\mathbf G=\prod_{i=1}^n\mathbf C_{p^{a_i}}$ has $p^n-1$ elements of order $p$, so it has $\frac{p^n-1}{p-1}$ cyclic subgroups of order $p$. This proves $(1)$. To prove $(2)$ notice that for $M\subseteq \mathbb N_n$ there are $(p-1)^{|M|}$ of elements of order $p$ whose $i$-th projection is non-trivial  if and only if  $i\in M$. This implies that there are $(p-1)^{|M|-1}$ cyclic subgroups of $\mathbf G$ of order $p$ whose $i$-th projection is non-trivial  if and only if  $i\in M$.

By a root of degree $p$ of $x$ we mean any element $y$ with $y^p=x$; we also say that $y$ is a $p$-th root of $x$. Notice that in a cyclic $p$-group any element has $p$ roots of degree $p$ if it is not a generator of the cyclic group, and no $p$-th roots otherwise. Now let $\mathbf H=\lb \overline x\rb=\lb(x_1,x_2,...,x_n)\rb$, and let $\mathbf H_i$ be $i$-th projection of $\mathbf H$. Then $\mathbf H_i=\lb x_i\rb$ for all $i\leq n$, and we have: $\mathbf H$ is a non-maximal cyclic subgroup of $\mathbf G$  if and only if  $\overline x$ has a $p$-th root in $\mathbf G$  if and only if  $x_i$ has a $p$-th root in $\mathbf C_{p^{a_i}}$ for all $i\leq n$  if and only if  $|H_i|< p^{a_i}$ for all $i\leq n$. This proves $(3)$.

%then $\mathbf H_i=\lb x_i\rb$ for all $i\in n$. If each $\mathbf H_i$ is not a maximal cyclic subgroup of $\mathbf G$, i.e. $|H_i|<p^{a_i}$, then each $x_i$ has its $p$-th root $y_i$. Then $\mathbf H<\lb\overline y\rb=\lb(y_0,y_1,...,y_{n-1})\rb$. On the other hand, if $\mathbf H$ is not a maximal cyclic subgroup of $\mathbf G$, then its generator $\overline x$ has a $p$-th root, which implies that each $x_i$ has a $p$-th root.

Let $\mathbf H=\lb \overline x\rb=\lb(x_1,x_2,...,x_n)\rb$ be a non-maximal cyclic subgroup of $\mathbf G$ of order $p^h>1$, and let $\mathbf H_i$ be the $i$-th projection of $\mathbf H$ of order $p^{h_i}$ for all $i\leq n$. Then $\mathbf H_i=\lb x_i\rb$ for all $i\leq n$ . Now for each $i\leq n$ the group $\lb x_i\rb$ is a non-maximal cyclic subgroup of $\mathbf C_{p^{a_i}}$, and $x_i$ has $p$ roots of degree $p$, so $\overline x$ has $p^n$ roots of degree $p$. As in each cyclic subgroup of $\mathbf G$ of order $p^{h+1}$ containing $\lb\overline x\rb$ as a subgroup has  $p$ roots of degree $p$, it follows that $\lb\overline x\rb$ is subgroup of $p^{n-1}$ cyclic subgroups of order $p^{h+1}$. This proves the first part of $(4)$. Let us prove the second part as well.

Let $L=\{i\in \mathbb N_n\mid h_i>0\}$, and let $L\subseteq M\subseteq \mathbb N_n$. Then the number or $p$-th roots of $\overline x$, whose $i$-th projection is non-trivial  if and only if  $i\in M$, is $p^{|L|}(p-1)^{|M|-|L|}$. Again, since each cyclic subgroup of order $p^{h+1}$ containing $\mathbf H$ as a subgroup contains $p$ roots of degree $p$ of $\overline x$, $\mathbf H$ is a subgroup of $p^{|L|-1}(p-1)^{|M|-|L|}$ cyclic subgroups of order $p^{h+1}$ whose $i$-th projection is non-trivial  if and only if  $i\in M$. This proves the second part of $(4)$, and the theorem has been proved.
\end{dokaz}\\

Note that there are non-abelian finite groups whose enhanced power graphs are isomorphic to a $p$-semitree. The group with representation $\mathbf G = \langle x, y \mid x^3 = y^3 = [x, y]^3 = 1\rangle$, mentioned in \cite{power graph}, of order $27$ and exponent $3$ has enhanced power graph isomorphic to the 3-semitree $S_3(1,1,1)$. Another example would be $\mathbf M_{16}$ with representation $\mathbf M_{16}=\langle a,x\mid a^8=x^2=1,xax^{-1}=a^5\rangle$ whose enhanced power graph is isomorphic to $S_2(3,1)$.

\begin{teo}
A graph $\Gamma$ is the enhanced power graph of a finite abelian group  if and only if  $\Gamma\cong \bigboxtimes_{i=1}^nS_{p_i}(\overline a_i)$ for some tuples $\overline a_i$ of positive integers and for some pairwise different primes $p_1,p_2,...,p_n$.
\end{teo}

\begin{dokaz}
The theorem follows from Theorem \ref{karakterizacija} and Lemma \ref{proizvodi}.
\end{dokaz}\\

%The following are Theorems 2.29, 4.12, and 5.39 in \cite{rotman}:
%\begin{enumerate}
%\item[$(a)$] Let $\mathbf G$ be a group with normal subgroups $\mathbf H$ and $\mathbf K$. If $H K= G$ and $H\cap K=1$, then $\mathbf G\cong \mathbf H\times\mathbf K$.
%\item[$(b)$] If $\mathbf P$ is a $p$-Sylow subgroup of a finite group $\mathbf G$, then all $p$-Sylow subgroups of $\mathbf G$ are conjugates of $\mathbf P$.
%\item[$(c)$] A finite group $\mathbf G$ is nilpotent iff it is the direct product of its Sylow subgroups.
%\end{enumerate}
%Then obviously, if $\mathbf G$ is nilpotent, then its $p$-Sylow subgroup is unique for each prime divisor $p$ of $|G|$. On the other hand, if $\mathbf G$ has a unique $p$-Sylow subgroup $\mathbf S_p$, then by $(b)$ $\mathbf S_p\trianglelefteq\mathbf G$. Since the same holds for all prime divisors of $|G|$, by $(a)$ $\mathbf G$ is direct product of its Sylow subgroups. This gives us the following theorem.

In the continuation of this section we will give a graph property which, for any finite group $\mathbf G$, is satisfied by $\mathcal G_e(\mathbf G)$  if and only if  $\mathbf G$ is nilpotent.

The following theorem is a well known fact about nilpotent groups.

\begin{teo}\label{nilpotentne grupe}
Let $\mathbf G$ be a finite group. Then the following conditions are equivalent:
\begin{enumerate}
\item $\mathbf G$ is nilpotent;
\item $\mathbf G$ is the direct product of its Sylow subgroups;
\item $\mathbf G$ has a unique $p$-Sylow subgroup, for all prime divisors $p$ of $|G|$.
\end{enumerate}
\end{teo}

%For proof of Theorem \ref{nilpotentne grupe} refer to Theorems 2.29, 4.12, and 5.39 in \cite{rotman}.
% Our goal now is, given the enhanced power graph of a group, to conclude as much as possible about the group itself. %Namely, Proposition \ref{prebrajanje} gives us a criterion for deciding whether the group in question is nilpotent. The next statement follows from Proposition \ref{prebrajanje} and Theorem \ref{nilpotentne grupe}.%we  know from the group's enhanced power graph whether the group is nilpotent or not. Namely, a group $\mathbf G$ of order $\prod_{i=0}^{n-1}p_i^{k_i}$ is nilpotent iff, for each $i\in n$, $\mathbf G$ has $p_i^{k_i}$ of elements whose order is a power of $p_i$.\footnote{Videti da li ovo da formuli\v semo kao posledicu nekako.}
In the next proposition we will prove that from a finite group's enhanced power graph we can conclude whether the group is nilpotent.

\begin{tvr}\label{sivljovesvost na osnovu brojeva}
Let $\mathbf G$ be a finite group of order $p_1^{k_1}p_2^{k_2}...p_n^{k_n}$, where $p_1,p_2,...,p_n$ are pairwise different primes. For each $i\leq n$ and $j\in\{0,1,...,k_i\}$, let $\mathcal D_i^j$ be the set of all maximal cliques of $\mathcal G_e(\mathbf G)$ of sizes divisible by $p_i^j$. Then:
\begin{enumerate}
\item $\mathbf G$ has a unique $p_i$-Sylow subgroup  if and only if  $\sum_{j=0}^{k_i}\big(\varphi(p_i^j)\cdot\lvert \mathcal D_i^j/\approx_{p_i^j}\rvert\big)=p_i^{k_i}$;
\item $\mathbf G$ is nilpotent  if and only if  $\sum_{j=0}^{k_i}\big(\varphi(p_i^j)\cdot|\mathcal D_i^j/\approx_{p_i^j}|\big)=p_i^{k_i}$ for all $i\leq n$. 
\end{enumerate}
\end{tvr}

\begin{dokaz}
Let $i\leq n$, and suppose that $\sum_{j=0}^{k_i}\big(\varphi(p_i^j)\cdot|\mathcal D_i^j/\approx_{p_i^j}|\big)=p_i^{k_i}$. By Proposition \ref{prebrajanje}, for any $j$, $\varphi(p_i^j)\cdot|\mathcal D_i^j/\approx_{p_i^j}|$ is the number of elements of order $p_i^j$ in $G$. Therefore, $\sum_{j=0}^{k_i}\big(\varphi(p_i^j)\cdot|\mathcal D_i^j/\approx_{p_i^j}|\big)$ is the number of elements in $G$ whose orders are powers of $p_i$. Since the number of such elements is $p_i^{k_i}$  if and only if  $\mathbf G$ has a unique $p_i$-Sylow subgroup, the first statement has been proved.

Now obviously, the equality holds for all $i\leq n$  if and only if  $\mathbf G$ has a unique $p_i$-Sylow subgroup for all $i\leq n$, which is, by Theorem \ref{nilpotentne grupe}, true if and only if  $\mathbf G$ is nilpotent.
\end{dokaz}\\

The following lemma provides us some necessary conditions for a finite graph to be the enhanced power graph of a group.

\begin{lem}\label{potrebni uslovi}
Let $\Gamma$ be a finite graph of order $p_1^{k_1}p_2^{k_2}...p_n^{k_n}$. Let $\mathcal Cl=\{C_1,C_2,...,C_n\}$ be the set of all maximal cliques of $\Gamma$, and let $\mathcal B=\{\bigcap_{i\in M}C_i\mid M\in\mathcal P^+(\mathbb N_n)\}$ be the set of all intersections of maximal cliques. If $\Gamma$ is the enhanced power graph of a finite group, then:
\begin{enumerate}
\myitem{$(E1)$} $\bigcap\{C\mid C\in\mathcal Cl\}\neq\emptyset$;\label{neprazan presek}
\myitem{$(E2)$} $\lvert C\rvert\big| p_1^{k_1}p_2^{k_2}...p_n^{k_n}$ for any $C\in \mathcal Cl$;\label{grupa deljiva klikom}
\myitem{$(E3)$} For any $B_1,B_2\in\mathcal B$: If $B_1,B_2\subseteq C$ for some $C\in\mathcal Cl$, then $|B_1\cap B_2|=\gcd(|B_1|,|B_2|)$.\label{presek je nzd}
\end{enumerate}
\end{lem}

\begin{dokaz}
Suppose that $\Gamma=\mathcal G_e(\mathbf G)$ for a finite group $\mathbf G$. By Lemma \ref{klike i ciklicne grupe} all elements of $\mathcal B$ and $\mathcal Cl$ are cyclic subgroups of $\mathbf G$. \ref{neprazan presek} follows from the fact that all subgroups of $\mathbf G$ contain the identity element of $\mathbf G$. $|G|$ is divisible by order of any subgroup of $\mathbf G$, which implies \ref{grupa deljiva klikom}.

Suppose that $B_1,B_2\in\mathcal B$ and $B_1,B_2\subseteq C$ for some $C\in\mathcal Cl$. Then $\mathbf B_1$, $\mathbf B_2$ and $\mathbf C$ are cyclic groups with $\mathbf B_1,\mathbf B_2\leq\mathbf C$, and the order of $\mathbf B_1\cap\mathbf B_2$ is $\gcd(|B_1|,|B_2|)$. This proves \ref{presek je nzd}.
%
%Then $\mathbf B_1\cap\mathbf B_2\leq \mathbf B_1,\mathbf B_2$, and therefore $|B_1\cap B_2|\ \big|\ \gcd(|B_1|,|B_2|)$. On the other hand, if $\mathbf D$ is the unique cyclic subgroup of $\mathbf C$ of order $\gcd(|B_1|,|B_2|)$, then $\mathbf D\leq\mathbf B_1,\mathbf B_2$, which implies $\mathbf D\leq \mathbf B_1\cap\mathbf B_2$. Therefore $\gcd(|B_1|,|B_2|)\ \big|\ |B_1\cap B_2|$, and thus $|B_1\cap B_2|=\gcd(|B_1|,|B_2|)$.
\end{dokaz}\\

Notice that for any finite graph, with the set of maximal cliques $\mathcal Cl$, and the set of all intersections of maximal cliques $\mathcal B$, \ref{presek je nzd} implies the following two conditions:
\begin{enumerate}
\myitem{$(E4)$} For any $B_1,B_2\in\mathcal B$: $B_1\subseteq B_2$ implies $|B_1|\big||B_2|$;\label{deljiv manjim}
\myitem{$(E5)$} For any $B_1,B_2\in\mathcal B$: If $|B_1|\big||B_2|$ and $B_1,B_2\subseteq C$ for some $C\in\mathcal Cl$, then $B_1\subseteq B_2$.\label{delilac je manji}
\end{enumerate}
Therefore, \ref{deljiv manjim} and \ref{delilac je manji} are also necessary conditions for a finite graph to be the enhanced power graph of a group.

\begin{figure}[h]
\begin{center}
\caption{}\label{nacrtan graf}
\includegraphics[width=70mm]{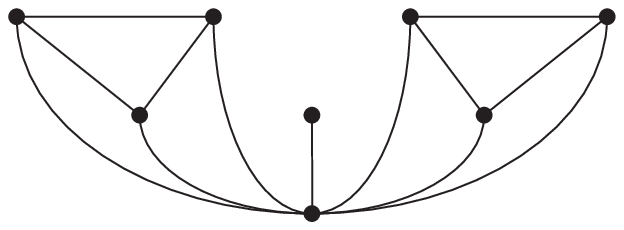}
\end{center}
\end{figure}

However, conditions \ref{neprazan presek}-\ref{presek je nzd} are not sufficient for a graph to be the enhanced power graph. For example, the graph in Figure \ref{nacrtan graf} satisfies \ref{neprazan presek}-\ref{presek je nzd}, although it is not enhanced power graph of any finite group. Namely, $\mathbf C_8$, $\mathbf C_4\times\mathbf C_2$, $\mathbf C_2\times\mathbf C_2\times\mathbf C_2$, $\mathbf D_8$, and $\mathbf Q_8$ are the only groups of order 8, and the graph from Figure  \ref{nacrtan graf} is isomorphic to none of their enhanced power graphs.

For a group $\mathbf G$ and a prime $p$ we shall denote with $G_p$ the set $\{x\in G\mid (\exists i\in\mathbb N_0)(o(x)=p^i)\}$. Note that $G_p$ may not be a subgroup of $\mathbf G$. However, if $\mathbf G$ has the unique $p$-Sylow subgroup $\mathbf P$, then $\mathbf P=\mathbf G_p$.

We will prove that the enhanced power graph of a finite group with a unique $p$-Sylow subgroup determines the enhanced power graph of its $p$-Sylow subgroup. For this reason we will introduce the notion of $p$-component of a graph.

\begin{de}
Let $\Gamma$ be a finite graph, and let $p$ be a prime number. Let $\mathcal Cl=\{C_1,C_2,...,C_n\}$ be the set of all maximal cliques of $\Gamma$, and let $\mathcal B=\{\bigcap_{i\in M}C_i\mid M\in\mathcal P^+(\mathbb N_n)\}$ be the set of all intersections of maximal cliques of $\Gamma$. A {\bf $p$-component} of $\Gamma$, denoted with $\Gamma_p$, is any induced subgraph of $\Gamma$ with property:
\begin{align}
\begin{split}
\text{For each $B\in\mathcal B$, of size $p^{k}l$ with $p\not|l$,}\\
\text{$B$ contains exactly $p^{k}$ vertices of $\Gamma_p$.}\label{svojstvo za silova}
\end{split}
\end{align}
\end{de}

Now we give an algorithm which for the enhanced power graph of a finite group having a unique $p$-Sylow subgroup returns the enhanced power graph of the $p$-Sylow subgroup.

\begin{tvr}\label{algoritam za p komponentu}
Let $\Gamma$ be a finite graph satisfying conditions \ref{neprazan presek} and \ref{presek je nzd}, and let $p$ be a prime divisor of  $|\Gamma|$. Then there is an algorithm for constructing a $p$-component of $\Gamma$.
\end{tvr}

\begin{dokaz}
Let $C_1,C_2,...,C_n$ be all maximal cliques of $\Gamma$, and let us denote $\bigcap_{i\in M}C_i$ with $B_M$ for every $M\in\mathcal P^+(\mathbb N_n)$. Let $n_{[p]}=\max\{p^k\mid k\in\mathbb N_0\text{ and }p^k|n\}$ for any $n\in\mathbb N$.  We claim that the following algorithm constructs a $p$-component of $\Gamma$.
\begin{enumerate}
\item Mark $|B_{\mathbb N_n}|_{[p]}$ vertices from $B_{\mathbb N_n}$, and cross out $\mathbb N_n$.
\item Choose a set $M_1\in\mathcal P^+(\mathbb N_n)$  whose all supersets but itself have been crossed out. Then mark more vertices from $B_{M_1}\setminus\bigcup_{M\supsetneq M_1}B_M$, if any is needed, in order for $B_{M_1}$ to contain exactly $|B_{M_1}|_{[p]}$ marked vertices, and cross out $M_1$.
%
%Choose a set $M_1\in\mathcal P^+(\mathbb N_n)$ for which there are less than $|B_{M_1}|_{[p]}$ marked vertices in $B_{M_1}$, and such that for all $M\supsetneq M_1$ there are exactly  $|B_M|_{[p]}$ marked vertices in $B_M$. Then mark more vertices from $B_{M_1}\setminus\bigcup_{M\supsetneq M_1}B_M$ in order for $B_{M_1}$ to contain exactly $|B_{M_1}|_{[p]}$ marked vertices.
\item Repeat the step 2. until all non-empty subsets of $\mathbb N_n$ have been crossed out.
%
%Repeat the step 2. until there are exactly $|B_M|_{[p]}$ marked elements in $B_M$ for all $M\in\mathcal P^+(\mathbb N_n)$.
\item Construct $\Gamma'$ as subgraph of $\Gamma$ induced by the set of all marked vertices.
%After we finish marking vertices from all intersections of maximal cliques in this manner, the subgraph $\Gamma_p$ induced by all marked vertices is a $p$-component of $\Gamma$.
\end{enumerate}

It is sufficient to prove that, for any $M_1\in\mathcal P^+$ whose all supersets but itself have been crossed out, $|B_{M_1}|_{[p]}$ is at least the sum of $|B_{M_1}\setminus\bigcup_{M\supsetneq M_1}B_i|$ and the number of vertices from $\bigcup_{M\supsetneq M_1}B_i$ that have been marked. Namely, were it true, then the algorithm finishes, and it outputs $\Gamma'$ which satisfies (\ref{svojstvo za silova}).

Let $M_1\in\mathcal P^+(\mathbb N_n)$, and let us denote all subsets $B_M$ for which $M_1\subsetneq M\subseteq \mathbb N_n$ with $B_1,B_2,...,B_m$. Suppose that all strict supersets of $M_1$ have been crossed out, but that $M_1$ has not. Then each $B_i$ contains exactly $|B_i|_{[p]}$ marked vertices. Let $p^{k_1}=|B_{M_1}|_{[p]}$, and let $p^{k_2}=\max\{|B_i|_{[p]}\ \big|\ i\leq m\}$. Notice that  $p^{k_2}\leq p^{k_1}$ by \ref{deljiv manjim}. Without loss of generality let $|B_1|_{[p]}=p^{k_2}$. We claim that all marked vertices from each $B_i$ belong to $B_1$. Namely, for all $i\leq m$ we have $B_1,B_i\subseteq B_{M_1}\subseteq C_j$ for some $j\leq n$, so by \ref{presek je nzd} it follows $|B_1\cap B_i|_{[p]}=|B_i|_{[p]}$. Therefore, all marked vertices from $B_i$ belong to $B_i\cap B_1\subseteq B_1$, which proves the claim. For this reason, $\bigcup_{i=1}^mB_i$ contains exactly $p^{k_2}$ marked vertices, and it remains for us to show that $p^{k_1}\leq |B_{M_1}\setminus\bigcup_{i=1}^mB_i|+p^{k_2}$.

Let $\mathbf G$ be a cyclic group of order $|B_{M_1}|$, and let $\mathbf G_i$ be the subgroup of $\mathbf G$ of order $|B_i|$ for all $i\leq m$. Notice that $|G_i\cap C_j|=\gcd(|G_i|,|G_j|)$. Then by \ref{presek je nzd} and by the inclusion-exclusion principle we get
\begin{align*}
&|B_1\cup B_2\cup ...\cup B_m|=\sum_{i=1}^m(-1)^{i+1}\sum_{J\in\mathbb N_n^{[i]}}\Big|\bigcap_{j\in J}B_i\Big|=\sum_{i=1}^m(-1)^{i+1}\sum_{J\in\mathbb N_n^{[i]}}\Big|\bigcap_{j\in J}G_i\Big|\\
&=\Big|\bigcup_{i=1}^nG_i\Big|\leq |G|-(p^{k_1}-p^{k_2})=|B_{M_1}|-(p^{k_1}-p^{k_2}),
\end{align*}
because $\bigcup_{i=1}^mG_i$ contains an element of $\mathbf G$ of order $p^i$ for no $i>k_2$. Therefore,  $p^{k_1}\leq |B_{M_1}\setminus\bigcup_{i=1}^mB_i|+p^{k_2}$, and it follows that one can mark $p^{k_1}-p^{k_2}$ vertices from $B_{M_1}\setminus\bigcup_{M\supsetneq M_1}B_M$, and then for any $M\supseteq M_1$ there would be exactly $|B_M|_{[p]}$ marked vertices in $B_M$.

Therefore, the algorithm eventually crosses over all non-empty subsets of $\mathbb N_n$, and it marks vertices of $\Gamma$ so that each intersection of maximal cliques $B$ contains exactly $|B|_{|p]}$ marked vertices. Then the subgraph $\Gamma'$ of $\Gamma$ induced by the set of marked elements is a $p$-component of $\Gamma$.
\end{dokaz}

\begin{lem}
Let $\mathbf G$ be a finite group, and let $p$ be a prime divisor of $|G|$. Then the subgraph of $\mathcal G_e(\mathbf G)$ induced by $G_p$ is a $p$-component of $\mathcal G_e(\mathbf G)$.
\end{lem}

\begin{dokaz}
Let $\mathcal Cl$ be the set of all maximal cliques of $\mathcal G_e(\mathbf G)$, and let $\mathcal B$ be the set of all intersections of maximal cliques of $\mathcal G_e(\mathbf G)$. By Lemma \ref{klike i ciklicne grupe} $\mathbf B$ is a cyclic subgroup of $\mathbf G$ for each $B\in\mathcal B$. Therefore, if $|B|=p^kl$ for some $k$ and $l$, $p{\not|}l$, then $\mathbf B$ has a unique subgroup of order $p^k$, so $B$ contains exactly $p^k$ elements whose orders are powers of $p$.
\end{dokaz}

\begin{lem}\label{za silova}
Let $\mathbf G$ be a finite group, let $p$ be a prime divisor of $|G|$, and let $\Gamma=\mathcal G_e(\mathbf G)$. %Let $\mathcal Cl=\{C_1,C_2,...,C_n\}$ be the set of all maximal cliques of $\Gamma$, and let $\mathcal B=\{\bigcap_{i\in M}C_i\mid M\in\mathcal P^+(\mathbb N_n)\}$ be the set of all intersections of maximal cliques of $\Gamma$.
 Then all $p$-components of $\Gamma$ are isomorphic to $\Gamma[G_p]$.
\end{lem}

\begin{dokaz}
Let $\Delta_1$ and $\Delta_2$ be two $p$-components of $\Gamma$. For any $x,y\in G$ we say $x\equive y$  if and only if  $\overline N_\Gamma(x)=\overline N_\Gamma(y)$. By Lemma \ref{uredjene klase} ({\it 2.}), $x\equive y$  if and only if  $x\in C\Leftrightarrow y\in C$ for all $C\in\mathcal Cl$. By Lemma \ref{uredjene klase} ({\it 4.}) $G/\equive$ can be ordered in the following way: $[x]_{\equive}\leq[y]_{\equive}$  if $y\in C\Rightarrow x\in C$ for all $C\in\mathcal Cl$. By Lemma \ref{uredjene klase} ({\it 1.} and {\it 5.})  $[x]_{\equive}\leq[y]_{\equive}$  if and only if  $\overline N_\Gamma(y)\subseteq \overline N_\Gamma(x)$  if and only if  $\langle [y]_{\equive}\rangle\leq\langle[x]_{\equive}\rangle$.

Height of an element $x$ of a poset is the maximal length of a chain whose greatest element is $x$.  By induction by height of $\equive$-class in the poset $(G/\equive,\leq)$, we prove that each $\equive$-class contains the same number of vertices of $\Delta_1$ and $\Delta_2$. The least element of $G/\equive$ is $\bigcap_{i=1}^nC_i$, and by (\ref{svojstvo za silova}) it has the same numbers of vertices of $\Delta_1$ and $\Delta_2$. Now let $D\in G/\equive$. Then by (\ref{svojstvo za silova}) $\langle D\rangle=\bigcap\{C\in\mathcal Cl\mid D\subseteq C\}$ has the same numbers of vertices of $\Delta_1$ and $\Delta_2$. By the induction hypothesis all $\equive$-classes $E<D$ have the same numbers of vertices of $\Delta_1$ and  $\Delta_2$. Note that for all such $E$ we have $E\subseteq\langle D\rangle$. This implies that $D$ contains the same numbers of vertices of $\Delta_1$ and $\Delta_2$. This proves the claim that each  $\equive$-class contains the same number of vertices of $\Delta_1$ and $\Delta_2$.

Now, because the transposition of $x$ and $y$ is an automorphism of $\Gamma$ whenever $x\equive y$, it follows that $\Delta_1$ and $\Delta_2$ are isomorphic by composition of transpositions of vertices belonging to same $\equive$-classes. Since $\Gamma(G_p)$ is a $p$-component of $\Gamma$, the lemma has been proven.
\end{dokaz}

\begin{teo}
Let $\mathbf G$ be a finite group with a unique $p$-Sylow subgroup $\mathbf G_p$ for a prime $p$.
%, and let $\Gamma=\mathcal G_e(\mathbf G)$.
 Then $\mathcal G_e(\mathbf G_p)$ is isomorphic to the $p$-component of $\mathcal G_e(\mathbf G)$.
%$\cong\Gamma_p$, where $\Gamma_p$ is the graph satisfying (\ref{svojstvo za silova}).
\end{teo}

\begin{dokaz}
The theorem follows directly from Lemma \ref{za silova}.
\end{dokaz}\\

If a finite graph $\Gamma$ of size $p_1^{n_1}p_2^{n_2}...p_k^{n_k}$ satisfies \ref{neprazan presek}-\ref{presek je nzd}, then a $p_i$-component $\Gamma_i$ can be constructed as described in Proposition \ref{algoritam za p komponentu} for any $i\leq k$. In this case $\Gamma$ is the enhanced power graph of a finite abelian group  if and only if  $\Gamma\cong\bigboxtimes_{i=1}^k\Gamma_{p_i}$ where $\Gamma_{p_i}$ is a $p_i$-semitrees for each $i$. If there are some $\Gamma_{p_i}$'s which are not $p$-semitrees, then we can not tell whether $\Gamma$ is an  enhanced power graph of a finite group.

%Let $p=p_i$ for some $i\leq n$. Let us see how to check whether $\Gamma_p$ is a $p$-semitree. Let $\mathcal Cl^p=\{C_1^p,C_2^p,...,C_m^p\}$ be the set of all maximal cliques of $\Gamma_p$, and let $\mathcal B^p=\{\bigcap_{i\in M}C_i^p\mid M\in\mathcal P^+(\mathbb N_m)\}$ be the set of all intersections of maximal cliques of $\Gamma$. Then $|B|$ is a power of the prime $p$ for any $B\in\mathcal B^p$, and every $B_1\in\mathcal B^p$ of size $p^{k+1}$ contains a set $B_2\in\mathcal B^p$ of size $p^k$. If $\Gamma_p$ was a $p$-semitree, i.e. an enhanced power graph of a finite abelian group, then 

To check whether $\Gamma_p$ is a $p$-semitree, i.e. the enhanced power graph of a finite abelian $p$-group, we can, as in Lemma \ref{prebrajanje}, figure how many elements of any order would the group have. From this information we can, as in Proposition  \ref{prepoznavanje na osnovu prebrajanja}, conclude to which tuple $\overline a$ would $\Gamma_p$ be related as a $p$-semitree. Then remains to check whether $\Gamma_p\cong S_p(\overline a)$.

If $\Gamma$ does not satisfy any of \ref{neprazan presek}-\ref{presek je nzd}, then it is not the enhanced power graph of any group.

%\begin{enumerate}
%\item For any $i\leq n$, there is the unique $\Gamma_{p_i}$, and $|\Gamma_{p_i}|\geq p_i^{k_i}$. 
%\end{enumerate}

\section{Perfectness of enhanced power graph}
\label{perfektnost}

 A graph $\Gamma$ is called perfect if for every finite induced subgraph $\Delta$ of $\Gamma$ the chromatic number of $\Delta$ is equal to the maximum size of a clique of $\Delta$. A graph $\Gamma$ is Berge graph if neither $\Gamma$ nor its complement contain an odd-length cycle of size at least 5 as an induced subgraph. By The Strong Perfect Graph Theorem%\footnote{Mogao bih da ubacim gde mo\v ze da se na\dj e dokaz za tu teoremu u literaturu.}
, which was proved in \cite{strong perfect graph}, a graph $\Gamma$ is perfect  if and only if  $\Gamma$ is a Berge graph. For that reason, the perfect graphs and the Berge graphs make up the same class of graphs. We shall denote the clique number of $\Gamma$ with $\omega(\Gamma)$, and the chromatic number of $\Gamma$ with $\chi(\Gamma)$. A graph is called weakly perfect if $\chi(\Gamma)=\omega(\Gamma)$.

\begin{pr}
There is an abelian finite group whose enhanced power graph is not perfect.
\end{pr}

\begin{dokaz}
Let $\mathbf G= \lb a_1\rb\times \lb a_2\rb\times \lb b_1\rb\times \lb b_2\rb\times \lb c_1\rb\times \lb c_2\rb$, where orders of $a_1$ and $a_2$, $b_1$ and $b_2$, and $c_1$ and $c_2$ are $2$, $3$, and $5$, respectively. Then its elements $a_1b_1$, $b_1c_2$, $a_2$, $b_2$, and $a_1c_1$ form the pentagon $\Pi$ as an induced subgraph of $\mathcal G_e(\mathbf G)$. Indeed, we have $a_1b_1,b_1c_2\in\lb a_1b_1c_2 \rb$, $b_1c_2,a_2\in\lb a_2b_1c_2\rb$, $a_2,b_2\in\lb a_2b_2\rb$, $b_2,a_1c_1\in \lb a_1b_2c_1\rb$, and $a_1c_1,a_1b_1\in\lb a_1b_1c_1\rb$, and no pair among $a_1b_1$ and $a_2$, $a_2$ and $a_1c_1$, $a_1c_1$ and $b_1c_2$, $b_1c_2$ and $b_2$, and $b_2$ and $a_1b_1$ is contained in a cyclic subgroup of $\mathbf G$. For instance, $\lb a_1c_1\rb$ and $\lb b_1c_2\rb$ contain two distinct cyclic subgroups of order $3$, $\lb c_1\rb$ and $\lb c_2\rb$, so any cyclic group containing $a_1c_1$ and $b_2c_2$ would contain two cyclic subgroups of the same order, which is impossible. Now maximal size of a clique in $\Pi$ is $2$, while chromatic number of $\Pi$ is $3$.
\end{dokaz}\ \\

Note that Theorem \ref{glavna za nilpotentne i perfektnost} shall imply that the above example gives a nilpotent group with non-perfect enhanced power graph of minimal size. Let us now give a complete characterization of finite nilpotent groups that have perfect enhanced power graph.

For a group $\mathbf G$, let us define relation $\simeqe_{\mathbf G}$ on $G$ so that $x\simeqe_{\mathbf G} y$ if $x=y$ or $x\sime y$.

%\begin{lem}\label{3d naocare}
%Let $\mathbf G=\mathbf H_1\times \mathbf  H_2\times ...\times \mathbf H_n$ be a finite group so that $(|H_i|,|H_j|)=1$ for all $i\neq j$. Then for any $\overline x=(x_1,x_2,...,x_n)$ and $\overline y=(y_1,y_2,...,y_n)$ from $G$ holds $\overline x\sim_{\mathbf G}\overline y$ iff $x_i\sim_{\mathbf H_i}y_i$ holds for all $i$.
%\end{lem}

%\begin{dokaz}
%If $\overline x\sim_{\mathbf G}\overline y$ holds, then there exists $\overline c\in G$ so that $\overline x,\overline y\in\lb \overline c\rb$. Obviously, if $\overline c=(c_1,c_2,...,c_n)$, then we have $x_i,y_i\in \lb c_i\rb$ for all $i$, which implies $x_i\sim_{\mathbf H_i}y_i$ for all $i$.

%If $x_i\sim_{\mathbf H_i}y_i$ holds for all $i$, then there exist $c_1\in H_1$, $c_2\in H_2$, ... ,$c_n\in H_n$ so that $x_i,y_i\in\lb c_i\rb$ for all $i$. Since orders of $\lb c_i\rb$ are pairwise relatively prime, then $\lb c_1\rb\times\lb c_2\rb\times ...\times\lb c_n\rb$ is a cyclic group containing both $\overline x$ and $\overline y$.
%\end{dokaz}

\begin{teo}\label{glavna za nilpotentne i perfektnost}
A finite nilpotent group has perfect enhanced power graph  if and only if  it has at most two non-cyclic Sylow subgroups.
\end{teo}

\begin{dokaz}
Let us assume that $\mathbf H_1$, $\mathbf H_2$ and $\mathbf H_3$ are three different non-cyclic Sylow subgroups of a nilpotent group $\mathbf G$. We shall prove that $\mathcal G_e(\mathbf H)$ contains a pentagon as induced subgraph, where $\mathbf H=\mathbf H_1\times \mathbf H_2\times \mathbf H_3$. By Theorem \ref{nilpotentne grupe} orders of $\mathbf H_1$, $\mathbf H_2$ and $\mathbf H_3$ are pairwise relatively prime. Let $a_1$ be an element of $H_1$ of maximal order, and let $a_2\in H_1\setminus\lb a_1\rb$. Then $a_1$ and $a_2$ are not adjacent in $\mathcal G_e(\mathbf H_1)$. In the same way can find elements $b_1,b_2\in H_2$, and $c_1,c_2\in H_3$ which are not adjacent in $\mathcal G_e(\mathbf H_2)$ and $\mathcal G_e(\mathbf H_3)$, respectively. By this, and by Lemma \ref{proizvodi}, it follows that the elements $a_1b_1$, $b_1c_2$, $a_2$, $b_2$, and $a_1c_1$ of $\mathbf H$ form a pentagon as induced subgraph of $\mathcal G_e(\mathbf H)$. This obviously implies that $\mathcal G_e(\mathbf G)$ also contains a pentagon as an induced subgraph. This proves one implication.

Suppose now that $\mathbf G=\prod_{i=1}^n\mathbf K_i$, where $\mathbf K_1,\mathbf K_2,...,\mathbf K_n$ have pairwise coprime orders. Assume that $\mathbf K_3,\mathbf K_4,...,\mathbf K_n$ are cyclic subgroups of $\mathbf G$, and let us denote $\mathbf K_1\times\mathbf K_2$ with $\mathbf K$. We shall prove that $\mathcal G_e(\mathbf K)$ does not contain any odd-length cycle of size at least $5$, nor a complement of an odd-length cycle of size at least $5$, as an induced subgraph. In the proof of this theorem we shall refer to these graphs as forbidden graphs.

First we shall note that $a\simeqe_{\mathbf P}b\simeqe_{\mathbf P}c\simeqe_{\mathbf P}d$, where $\mathbf P$ is finite a $p$-group $\mathbf P$, implies $a\simeqe_{\mathbf P}c$ or $b\simeqe_{\mathbf P}d$. Indeed, since $b$ and $c$ are contained in a same cyclic $p$-group, one is a power of the other. If $c$ is a power of $b$, then $a\simeqe_{\mathbf P}c$, while if $b$ is a power of $c$, then we have $b\simeqe_{\mathbf P}d$.

Suppose now that $\mathcal G_e(\mathbf K)$ contains an odd-length cycle of size at least 5 as an induced subgraph. Let the cycle be $a_1b_1, a_2b_2,...,a_{2k+1}b_{2k+1}$, where $a_i\in H_1$ and $b_i\in H_2$ for all $i$. Then we have
$$a_1b_1\sime_{\mathbf K}a_2b_2\sime_{\mathbf K}a_3b_3\sime_{\mathbf K}...\sime_{\mathbf K}a_{2k+1}b_{2k+1}\sime_{\mathbf K}a_1b_1,$$
which  by Lemma \ref{proizvodi} implies
\begin{align*}
a_1\simeqe_{\mathbf K_1}a_2\simeqe_{\mathbf K_1}a_3\simeqe_{\mathbf K_1}...\simeqe_{\mathbf K_1}a_{2k+1}\simeqe_{\mathbf K_1}a_1,&\text{ and}\\
b_1\simeqe_{\mathbf K_2}b_2\simeqe_{\mathbf K_2}b_3\simeqe_{\mathbf K_2}...\simeqe_{\mathbf K_2}b_{2k+1}\simeqe_{\mathbf K_2}b_1.&
\end{align*}
Also, we have $a_1b_1\not\simeqe_{\mathbf K}a_3b_3$, which implies $a_1\not\simeqe_{\mathbf K_1}a_3$ or $b_1\not\simeqe_{\mathbf K_2}b_3$. This way we get:
\begin{align*}
a_1\not\simeqe_{\mathbf K_1}a_3&\text{ or }b_1\not\simeqe_{\mathbf K_2}b_3,\\
a_2\not\simeqe_{\mathbf K_1}a_4&\text{ or }b_2\not\simeqe_{\mathbf K_2}b_4,\\
&...\\
a_{n-1}\not\simeqe_{\mathbf K_1}a_1&\text{ or }b_{n-1}\not\simeqe_{\mathbf K_2}b_1,\\
a_{2k+1}\not\simeqe_{\mathbf K_1}a_2&\text{ or }b_{2k+1}\not\simeqe_{\mathbf K_2}b_2,
\end{align*}
which can be observed as an odd-sized cycle colored into two colors: the cycle is $(1,3),(2,4),(3,5),...,(2k+1,2)$, while a pair $(i,i+2)$ is colored with $a$ only if $a_i\not\simeqe a_{i+_{2k+1}2}$, and it is colored with $b$ only if $b_i\not\simeqe b_{i+_{2k+1}2}$. As such cycle has two subsequent elements, without loss of generality we can assume that we have $a_1\not\simeqe_{\mathbf K_1}a_3$ and $a_2\not\simeqe_{\mathbf K_1}a_4$. On the other hand, as seen above, $a_1\simeqe_{\mathbf K_1}a_2\simeqe_{\mathbf K_1}a_3\simeqe_{\mathbf K_1}a_4$ implies $a_1\simeqe_{\mathbf K_1}a_3$ or $a_2\simeqe_{\mathbf K_1}a_4$. We've got a contradiction, so $\mathcal G_e(\mathbf K)$ does not contain an odd-length cycle of size at least five.

Next we show that $\mathcal G_e(\mathbf K)$ does not contain the complement of an odd-sized cycle of size at least 7. Note that we already have the result for the size 5, because complement of a cycle of size 5 is again a cycle. Assume that it does contain the complement of the cycle $a_1b_1, a_2b_2,...,a_{2k+1}b_{2k+1}$, where $a_i\in H_1$ and $b_i\in H_2$ for all $i$. Now we have:
\begin{align*}
a_1\not\simeqe_{\mathbf K_1}a_2&\text{ or }b_1\not\simeqe_{\mathbf K_2}b_2,\\
%a_3\not\simeqe_{\mathbf K_1}a_4&\text{ or }b_3\not\simeqe_{\mathbf K_2}b_4,\\
%&...\\
a_2\not\simeqe_{\mathbf K_1}a_3&\text{ or }b_2\not\simeqe_{\mathbf K_2}b_3,\\
&...\\
a_{2k+1}\not\simeqe_{\mathbf K_1}a_1&\text{ or }b_{2k+1}\not\simeqe_{\mathbf K_2}b_1,
\end{align*}
which can also be observed as a cycle of size at least 7 colored into two colors. Such a cycle has three pairwise non-adjacent elements, so there are two non-adjacent elements in the cycle colored with the same color. Therefore, without the loss of generality we can assume that $a_1\not\simeqe_{\mathbf K_1}a_2$ and $a_3\not\simeqe_{\mathbf K_1}a_4$. On the other hand, $a_2\simeqe_{\mathbf K_1}a_4\simeqe_{\mathbf K_1}a_1\simeqe_{\mathbf K_1}a_3$, and it implies $a_1\simeqe_{\mathbf K_1}a_2$ or $a_3\simeqe_{\mathbf K_1}a_4$. Again a contradiction. Therefore, $\mathcal G_e(\mathbf K)$ does not contain any of the forbidden subgraphs as an induced subgraph.

Notice that $\mathbf C=\prod_{i=3}^n\mathbf K_i$ is a cyclic group. Suppose that $\mathbf G=\mathbf K\times\mathbf C$ contains a forbidden induced subgraph. Now by Lemma \ref{proizvodi}, because $x\simeqe_{\mathbf C}y$ for any $x,y\in C$, it follows that $\mathbf K$ contains a forbidden graph, which is a contradiction. Therefore, no enhanced power graph of a nilpotent group with at most two non-cyclic Sylow subgroups contains a forbidden graph as an induced subgraph. So by The Strong Perfect Graph Theorem we get that any finite nilpotent group with at most two non-cyclic subgroups has a perfect enhanced power graph.
\end{dokaz}\\

On the other hand, the following theorem implies that there are many groups which do not have perfect enhanced power graph, but whose clique number and chromatic number are equal.

\begin{teo}\label{za sire od nilpotentnih}
If a finite group $\mathbf G$ has an element whose order is equal to the exponent of $\mathbf G$, then $\mathcal G_e(\mathbf G)$ is weakly perfect.
\end{teo}

\begin{dokaz}
Let $\mathbf G$ be a finite group which has an element whose order is equal to the exponent of $\mathbf G$. Let $x$ be such an element, and let $o(x)=n$. Then $\omega(\mathcal G_e(\mathbf G))=n$. We can color the elements of $\lb x\rb$ with $n$ colors. Now, for any  $\approx$-class $D$ of $G$ there is a $\approx$-class $C\subseteq\lb x\rb$, such that $C$ and $D$ contain elements of the same order. Therefore, we can color the elements of $D$ with the same colors the elements of $C$ were colored with. Suppose now that all edges of $\mathcal G_e(\mathbf G)$ have been colored in this manner, and suppose that $y$ and $z$ are colored with a same color. Then $o(y)=o(z)$ and $y\not\approx z$, which implies $y\not\sime z$. For that reason, in this coloring of the graph $\mathcal G_e(\mathbf G)$ there is no pair of adjacent vertices colored with a same color, and we have $\chi(\mathcal G_e(\mathbf G))=\omega(\mathcal G_e(\mathbf G))$.
\end{dokaz}

\begin{pos}
If $\mathbf G$ is a finite nilpotent group, then $\mathcal G_e(\mathbf G)$ is weakly perfect.
\end{pos}

\begin{dokaz}
Let $\mathbf G$ be a finite nilpotent group. Then it is a direct product of its Sylow subgroups, implying $\mathbf G$ has an element whose order is equal to the exponent of $\mathbf G$. Therefore, by Theorem \ref{za sire od nilpotentnih} it follows that $\mathcal G_e(\mathbf G)$ is weakly perfect.
\end{dokaz}

\noindent
Samir Zahirovi\' c

Department of Mathematics and Informatics, University of Novi Sad, Serbia

{\it e-mail}: samir.zahirovic@dmi.uns.ac.rs \\

\noindent
Ivica Bo\v snjak

Department of Mathematics and Informatics, University of Novi Sad, Serbia

{\it e-mail}: ivb@dmi.uns.ac.rs\\

\noindent
Roz\' alia Madar\' asz

Department of Mathematics and Informatics, University of Novi Sad, Serbia

{\it e-mail}: rozi@dmi.uns.ac.rs

\end{document}